\documentclass[12pt,a4paper,american]{article}
\usepackage[T1]{fontenc}
\usepackage{amstext}
\usepackage[authoryear]{natbib}
\usepackage{amsfonts}
\usepackage[fleqn]{amsmath}
\usepackage{nicefrac}
\usepackage{graphicx}
\usepackage{mathrsfs}
\usepackage{listings}
\usepackage{cases}
\usepackage{tabularx, booktabs}
\usepackage{setspace}
\usepackage{multirow}
\usepackage{subcaption}
\usepackage[toc,page]{appendix}

\usepackage{cleveref}

\usepackage{vmargin}
\setmarginsrb{2cm}{2cm}{2cm}{2cm}{1cm}{1cm}{1cm}{1cm}

\usepackage{psboxit}

\usepackage{rotfloat}
\floatstyle{ruled}
\newfloat{algorithm}{tbp}{loa}
\floatname{algorithm}{Algorithm}
\usepackage[fleqn]{amsmath}
\makeatletter
\usepackage{pstricks}
\usepackage{pst-node}
\makeatother

\begin{document}
	\title{Adaptive Partitioning Strategy for High-Dimensional Discrete Simulation-based Optimization Problems}
	\author{{\bf Jing Lu}\\
		Operations Research Center\\
		Massachusetts Institute of Technology, USA
		\and
		{\bf Tianli Zhou}\\
		Department of Civil and Environmental Engineering\\
		Massachusetts Institute of Technology, USA
		\and
		{\bf Carolina Osorio}\\
		Department of Decision Sciences\\
		HEC Montréal, Canada}
	\date{\today}
	\maketitle
	\section* {Abstract}
	In this paper, we introduce a technique to enhance the computational efficiency of solution algorithms for high-dimensional discrete simulation-based optimization problems. The technique is based on innovative adaptive partitioning strategies that partition the feasible region using solutions that has already been simulated as well as prior knowledge of the problem of interesting.
	We integrate the proposed strategies with the Empirical Stochastic Branch-and-Bound framework proposed by \cite{xu2013empirical}. This combination leads to a general-purpose discrete simulation-based optimization algorithm that is both globally convergent and has good small sample (finite-time) performance. 
	The proposed general-purpose discrete simulation-based optimization algorithm is validated on a synthetic discrete simulation-based optimization problem and is then used to address a real-world car-sharing fleet assignment problem. Experiment results show that the proposed strategy can increase the algorithm efficiency significantly. 

	\section{Introduction}
	%
	Simulation-based optimization (SO), also referred as optimization via simulation (OvS), is the optimization of the performance of a stochastic system, where the objective function and/or constraints can only be estimated through stochastic simulations. 
	SO is often used to study complex systems, such as transportation \citep{osorio2013simulation, zhouOso2019, lu2020probabilistic}, petroleum and gas production \citep{wang2013mixed, li2016mixed}, building performance analysis \citep{nguyen2014review}, semiconductor manufacturing \citep{lin2015simulation} and supply chain design \citep{xu2009flexibility}.

	Based on the problem feasible region structures, \citet{hong2009brief} divide SO problems into three categories: ranking and selection, continuous SO, and discrete SO. 
	A ranking and selection problem usually have a relatively small number of feasible solutions, and all feasible solutions can be simulated at least once and the (near) best solution is chosen with a given confidence interval (e.g., \citet{kim2006selecting}). 
	For a continuous SO problem, the decision variables are continuous, while for a discrete SO problem is one with discrete decision variables. Methods developed to tackle continuous SO problems include stochastic approximation methods (e.g., \citet{robbins1951stochastic,bhatnagar2011stochastic}) and surrogate-based methods where one or more analytical models are used to approximate a simulated value (e.g., \citet{barton2006metamodel, angun2009response, regis2013combining, wang2018constrained}). 
	\citet{osorio2013simulation} propose metamodel SO, which use a surrogate model that combine a problem specific analytical model with a general purpose surrogate model.

	For a discrete SO problem, the number of feasible solutions is usually too large for each one to be simulated. The existing algorithms usually seek to leverage anticipated spatial structure in the solution space such as the clustering of good feasible solutions \citep{nelson2010optimization}.
	Some methods that aim at identifying solutions with good performances within small sampling budgets do not guarantee global convergence. Examples of this type of algorithms include Heuristic Constrained Genetic Algorithm (HCGA) \citep{tsai2014genetic}, Convergent Optimization via Most-Promising-Area Stochastic Search (COMPASS) \citep{hong2006discrete}, Adaptive Hyperbox Algorithm (AHA) \citep{xu2013adaptive}, R-SPLINE \citep{wang2013integer} and cgR-SPLINE \citep{nagaraj2014stochastically}.
	Other discrete SO algorithms focus on global convergence to the optimal solution(s) with unlimited computational resources.
	\citet{shi2000nested, andradottir2006overview, xu2013empirical} are example methods that rely on iteratively partition the feasible region into subregions and mainly search within some promising subregions.
	Some methods use a probabilistic model to compute the probability to sample each feasible solution. Commonly used probabilistic model includes Gaussian process \citep{sun2014balancing} or Gaussian Markov random fields \citep{salemi2019gaussian}.
	Surrogate-based methods are also used in discrete SO. Commonly used surrogate models include low-order polynomial model, Kriging model \citep{kleijnen2009kriging} and radial-basis model \citep{wild2008orbit}.
	%
	When certain problem specific information is available, metamodel SO algorithms can be used to efficiently solve a discrete SO problem \citep{zhouOso2019, zhou2020}.
	Nevertheless, developing discrete SO algorithms that can efficiently tackle high-dimensional problems remains a challenge.
	%
	Detailed overviews of SO literature are provided by \citet{fu2002optimization, amaran2016simulation, bhosekar2018advances}. 
	
	\cite{xu2013empirical} propose an Empirical Stochastic Branch-and-Bound (ESB\&B) framework for discrete SO problems based on the nested partitions method of \citet{shi2000nested} and stochastic branch-and-bound \citep{norkin1998branch}.
	It takes advantage of the partitioning structure of stochastic branch-and-bound method and  empirically estimates the bounds based on sampled solutions. The ESB\&B algorithm also uses improvement bounds to represent the potential of each subregion to guide the sampling strategy in the next iteration. The
	ESB\&B algorithm is globally convergent, i.e., given infinite simulation budget, it converges asymptotically to the global optimum. 
	As mentioned in \cite{xu2013empirical}, there are many valid partitioning strategies; however, a good partitioning strategy can usually improve the algorithm efficiency. Most partitioning strategies developed in the literature are generic and heuristic, for example, dividing the feasible region equally into $k$ subregions along a randomly chosen dimension \citep{xu2013empirical}. 
	
	In this paper, we propose an innovative adaptive partitioning strategy. It is embedded within the ESB\&B framework \citep{xu2013empirical} and forms a globally convergent algorithm for discrete SO problems. 
	The proposed partitioning strategy iteratively divides the feasible region in a fashion such that previously sampled solutions with similar performances are located in the same subregion. 
	It is an adaptive sample-based partitioning strategy that enhances small sample (finite-time) performances of the ESB\&B framework of \citet{xu2013empirical}. 
	This proposed partitioning strategy can take on problem-specific structures known a priori, such as spatial dependence in car-sharing fleet assignment problems, to further enhance the algorithm efficiencies without significant modifications. 
	
	This paper is organized as follows. Section \ref{sec:esbb} reviews the ESB\&B algorithm by \cite{xu2013empirical}. Section \ref{sec:optimalParition} presents the proposed adaptive partitioning strategy and its solving method developed by \citet{dunn2018optimal}. Section \ref{sec:evaluation} validates the proposed partitioning strategy on one synthetic and one real-world discrete SO problems. Section \ref{sec:conclude} concludes this paper.
	
	\section{ESB\&B framework} \label{sec:esbb}
	Let us first define the discrete SO problem.
	The goal is to find $\mathbf{x}$ that solves
	\begin{equation}
	\max_{\mathbf{x}\in \mathcal{X}} E[Y(\mathbf{x})]
	\end{equation}
	where $\mathbf{x} = [x_1,...,x_p]$ are the discrete decision variables, and  $\mathcal{X}$ is a \textbf{convex} feasible region that contains \textbf{finite} but a large number of feasible solutions, which can be represented by constraints of the form:
	\begin{align}
	&l_i\leq x_i\leq u_i, \quad i = 1,...,p\\
	&g_j(\mathbf{x})\leq 0, \quad j=1,...,q\\
	&l_i,x_i,u_i\in\mathbb{Z}, \quad i = 1,...,p.
	\end{align}
	The objective function $E[Y(\mathbf{x})]$ is the expected performance at point $\mathbf{x}$, which can only be estimated by generating observations of $Y(\mathbf{x})$ via simulation. 
	
	To solve the discrete SO problems of the above format, \citet{xu2013empirical} propose the Empirical Stochastic Branch-and-Bound (ESB\&B) framework, which converges to the globally optimal solution(s) asymptotically (i.e., under unlimited sampling budgets and simulation efforts). The ESB\&B algorithm is detailed in Algorithm \ref{alg:ESBB}.
	The algorithm terminates when the total sampling budget is used up. Whenever the algorithm is terminated, the final solution is the one with the maximum cumulative sample average. 
	
	\begin{algorithm}		
		\caption{ESB\&B framework \label{alg:ESBB}}
		\small
		\begin{enumerate}
			\item[0.] Initialization: set iteration counter $k = 0$, initial sampled solution set $\Phi^{0}=\emptyset$, initial partition $\mathcal{P}_0 = \{\mathcal{X}\}$, initial best subregion $\mathcal{R}^0 = \mathcal{X}$
			\item[1.] Partitioning:\\ If the best subregion $\mathcal{R}^k$ is singleton:
			\begin{enumerate}
				\item set the new full partition $\mathcal{P}^{'}_k = \mathcal{P}_k$			
			\end{enumerate}
			else:
			\begin{enumerate}
				\item construct a partition of the best subregion $\mathcal{P}(\mathcal{R}^k)$ 
				\item define the new full partition by $\mathcal{P}^{'}_k = (\mathcal{P}_k\textbackslash\{\mathcal{R}^k\})\cup\mathcal{P}(\mathcal{R}^k)$
				\item denote $\mathcal{X}^P$ the elements of $\mathcal{P}^{'}_k$ 
			\end{enumerate}
			\item[2.] Sampling and bounding:
			\begin{enumerate}
				\item[2.1] Solution sampling: 
				\begin{enumerate}
					\item for each subregion $\mathcal{X}^P\in\mathcal{P}(\mathcal{R}^k)$, randomly sample $\nu_R$ solutions
					\item if $k > 0$, for each subregion $\mathcal{X}^P\in \mathcal{P}_k\textbackslash\{\mathcal{R}^k\}$, sample $\theta(\mathcal{X}^P)$ solutions, where $\theta(\mathcal{X}^P)$ is computed in Step 2.3 of iteration $k-1$ based on information in $\Phi^{k-1}$
					\item aggregate all of the sampled solutions into set $S^k$, and set $\Phi^k = \Phi^{k-1}\cup S^k$		
				\end{enumerate} 
				\item[2.2] Bound estimation:
				\begin{enumerate}
					\item for each $\mathbf{x}\in S^k$, simulate $\Delta n_F$ replications if $\mathbf{x}\notin \Phi^{k-1}$, simulate $\Delta n_A$ additional replications if $\mathbf{x}\in \Phi^{k-1}$
					\item for $\mathcal{X}^P\in \mathcal{P}^{'}_k$, calculate estimated upper bound $\eta^{k+1}(\mathcal{X}^P)$
				\end{enumerate}
				\item[2.3] Sample allocation: compute the number of solutions to be sampled, $\theta(\mathcal{X}^P)$, for all $\mathcal{X}^P\in \mathcal{P}^{'}_k$ for iteration $k+1$ based on information in $\Phi^k$
			\end{enumerate}
			\item[3.] Updating partition and best subregion:
			\begin{enumerate}
				\item update the best subregion $\mathcal{R}^{k+1} = \arg\max_{\mathcal{X}^P\in\mathcal{P}^{'}_k}\{\eta^{k+1}(\mathcal{X}^P)\}$
				\item partition $\mathcal{P}_{k+1} = \mathcal{P}^{'}_k$
				\item set $k=k+1$, go to Step 1.
			\end{enumerate}
		\end{enumerate}
	\end{algorithm}
	
	There are two important steps within each iteration: i) sampling and bounding (Step 2), ii) and partitioning (Step 3). For details regarding the sampling and bounding step, we refer to Section 3 of \citet{xu2013empirical}. We discuss in details the partitioning step.

		\subsection*{Partitioning}
		Partitioning is the step that divides the estimated best subregion into a set of smaller subregions that are disjoint and nonempty. In the ESB\&B implementation of \citet{xu2013empirical}, the embedded generic partitioning strategy chooses (either deterministically or randomly) a variable $x_i$ and divide the best subregion equally along this dimension into $\omega$ disjoint subregions. A detailed description of this partitioning strategy is given by the online Appendix A of \citet{xu2013empirical}. 
		
		However, there exists many valid partitions for a given subregion.
		A good partitioning strategy can help locate the most promising subregion more efficiently, and hence allocate the sampling budget more efficiently. Let us take the following illustrative example to discuss the pros and cons of the current generic partitioning strategy used in ESB\&B and potential directions of improvement.
		
		\subsubsection*{Illustrative example}
\begin{figure}[htp]
		\minipage{1\textwidth}
		\centering
		\includegraphics[width=0.9\linewidth]{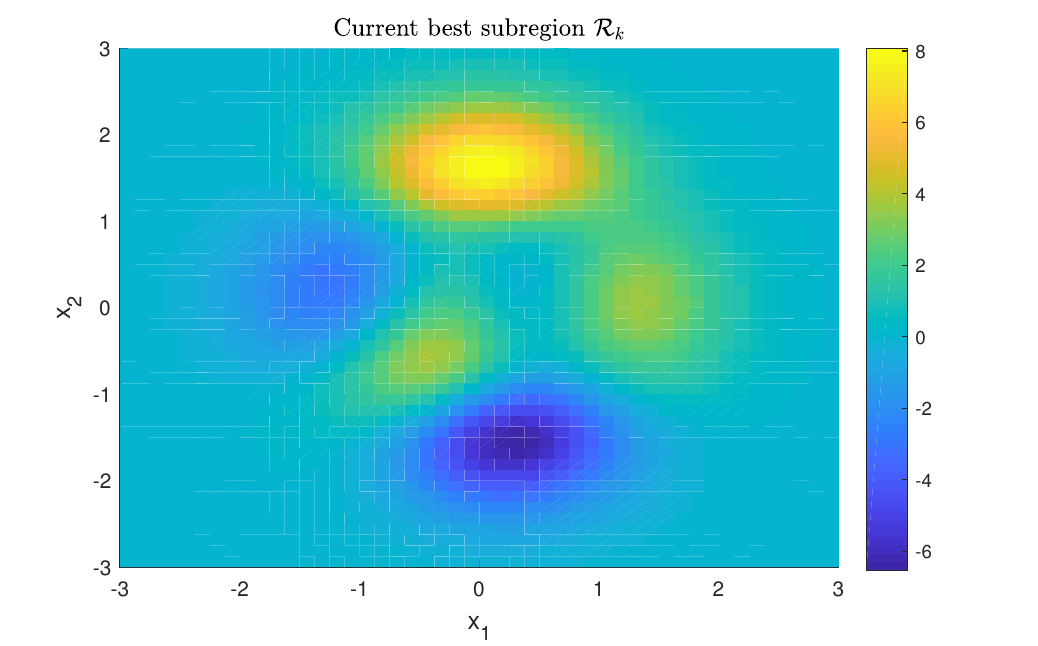}
		\caption{The ground truth values of $f(x_1,x_2)$ in the current best subregion.}
		\label{fig:chap4_1}
		\endminipage\hfill
	\end{figure}

\begin{figure}[htp]
		\minipage{1\textwidth}
		\centering
		\includegraphics[width=0.9\linewidth]{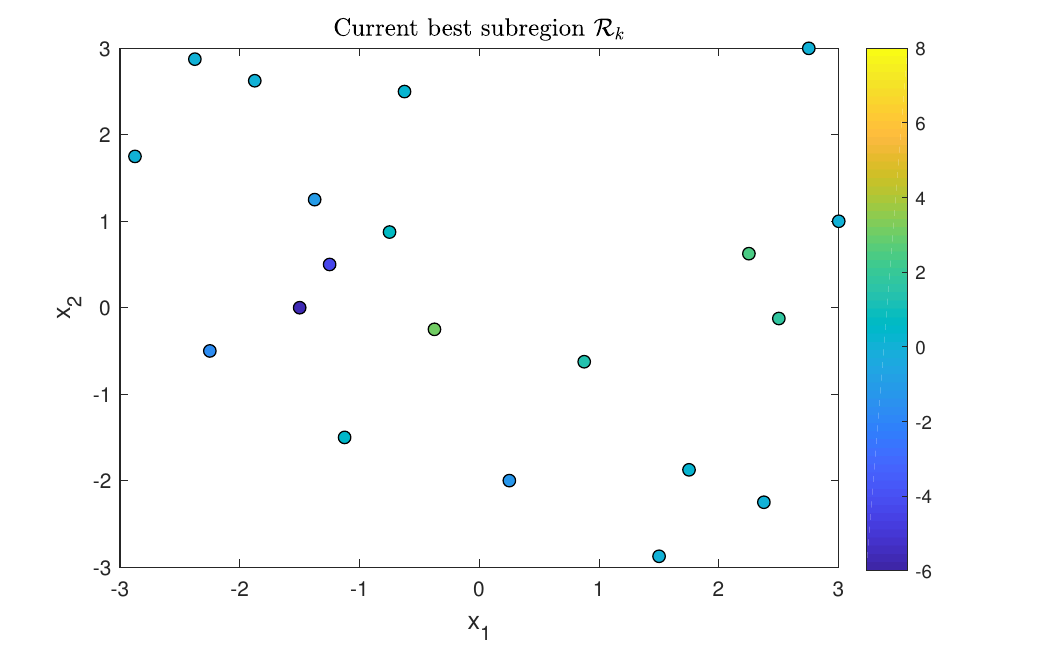}
		\caption{The sampled solutions of $f(x_1,x_2)$ in the current best subregion.}
		\label{fig:chap4_2}
		\endminipage\hfill
	\end{figure}

		Consider the maximization of a 2-dimensional deterministic function $f(x_1,x_2)$, for which neither the closed-form nor any structural information about the function is known. The ground truth of $f(x_1,x_2)$ in the current best subregion is shown in Figure \ref{fig:chap4_1}. However, the values of the function can only be known through sampling.
		Figure \ref{fig:chap4_2} shows samples that have already been evaluated in the current best subregion. 
		
\begin{figure}[htp]
		\minipage{1\textwidth}
		\includegraphics[width=1\textwidth]{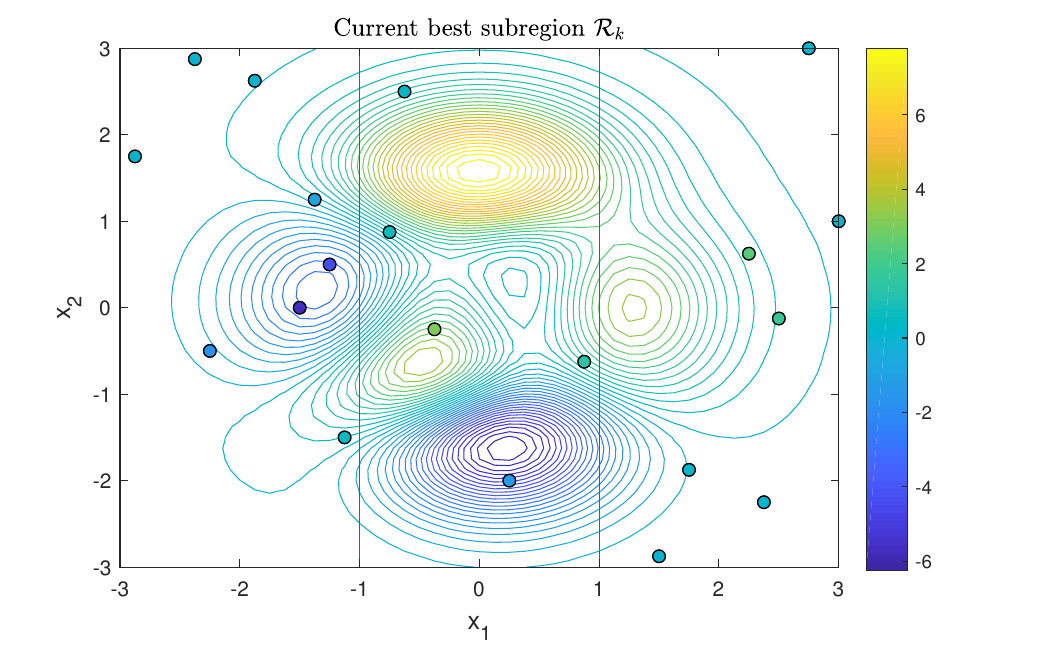}
		\subcaption{A naive partition of the current best subregion.}\label{fig:chap4_3}
		\endminipage\hfill
		\minipage{1\textwidth}
		\includegraphics[width=1\textwidth]{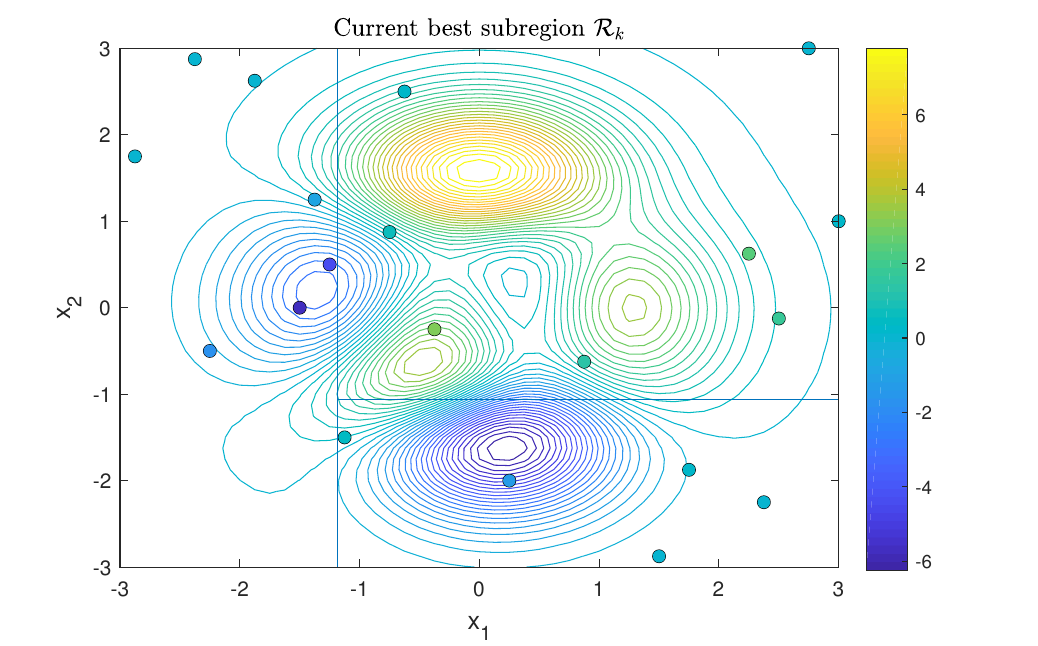}
		\subcaption{A better partition of the current best subregion.}\label{fig:chap4_4}
		\endminipage\hfill
		\caption{\label{fig:chap4partitioning} Different partitions of the current best subregion.}
	\end{figure}

		We need to further partition this region and sample from each new subregion. 
		Figure \ref{fig:chap4_3} shows the generic partitioning strategy used in the ESB\&B framework, which divides this region into equal parts along the chosen $x_1$ dimension. This partitioning strategy requires the users to pre-define the number of subregions to be divided into (denoted $\omega$), and it does not utilize any information from the solutions that have already been sampled.
		After partitioning, each new subregion gets an equal amount of sampling budget. Based on the performances of the sampled solutions, a new best subregion will be selected and further explored. In this case, the middle and right subregions both have the chance to be selected as the next best subregion as they both contains the peaks of the underlying function. If the right subregion is chosen as the next best subregion, it can take a while for the algorithm to return to explore the middle subregion where the true global maximum solution locates. This is because the sampling budget assigned to the best subregions are usually higher than that to the other subregions.
		
		On the other hand, there is another potential partition of the current best subregion given in Figure \ref{fig:chap4_4}. Note that this partition divides the subregion into three parts: two that contains the underlying function's basins and one contains all the peaks. This partition is better in the sense that it successfully identifies the patterns of the underlying function. It is obtained by grouping sampled solutions with similar performances in the same subregion. 
	
	In Section \ref{sec:optimalParition}, we propose an adaptive partitioning strategy which is defined by the previously sampled solutions in the subregion.
	
	\section{Adaptive partitioning strategy}\label{sec:optimalParition}
	The main idea underlying the proposed partitioning strategy is to find a partition of the current best subregion $\mathcal{R}^k$ in which the sampled solutions with similar performances are divided in the same subregion. 
	In this section, we introduce two sets of adaptive partitioning strategies: (i) parallel partition in Section \ref{sec:parallelPartition}, which applies to any problems, and (ii) hyperplane partition in Section \ref{sec:hyperplanePartition}, which applies when splitting features in the form of linear combinations of decision variables can be obtained from prior knowledge (e.g., spatial dependence in car-sharing fleet assignment problem).
	
	\subsection{Parallel partition}\label{sec:parallelPartition}
	A parallel partition is one that only contains cuts of the form $\mathbf{a}^T\mathbf{x}<b$ or $\mathbf{a}^T\mathbf{x}\geq b$ where $a_i\in\{0,1\}$ and $\sum_{i=1}^{p}a_i = 1$.
	The subregions resulting from a parallel partition are in the form of the intersection of a hyperbox with the best subregion. This family of partitions are favorable in the sense that it directly maps to the bounds of each decision variable and hence is easy to interpret and draw feasible uniform random samples from. A parallel partitioning strategy can be applied to any discrete SO problem, thus it is a generic partitioning strategy.
	
	We formulate the search for such a parallel partition as a minimization problem: 
	\begin{alignat}{3}
	\min_{\mathcal{P}(\mathcal{R}^k)} &\quad&  \sum_{\mathcal{X}^P\in\mathcal{P}(\mathcal{R}^k)}\sum_{\mathbf{x}\in \mathcal{X}^P\cap\Phi^{k}}(\bar{Y}(\mathbf{x})-\bar{Y})^2 & \label{eq:min}\\
	\text{s.t. } &\quad&  |\mathcal{X}^P\cap\Phi^{k}|\geq N_{min}, & \quad \forall \quad \mathcal{X}^P\in \mathcal{P}(\mathcal{R}^k), \label{eq:constrA}\\
	&\quad& |\mathcal{P}(\mathcal{R}^k)|\leq d, \label{eq:constrB}
	\end{alignat}
	where $\mathcal{P}(\mathcal{R}^k)$ is a valid parallel partition of the current best subregion $\mathcal{R}^k$, $\bar{Y}(\mathbf{x})$ is the cumulative sample mean of all $n(\mathbf{x})$ observations at solution $\mathbf{x}$:
	\begin{equation}
	\bar{Y}(\mathbf{x})= \frac{1}{n(\mathbf{x})}\sum_{s=1}^{n(\mathbf{x})}Y_s(\mathbf{x}).
	\end{equation}
	$\bar{Y}$ is average estimated mean performances of all sampled points in subregion $\mathcal{X}^P$:
	\begin{equation}\label{eq:partition_mean}
	\bar{Y} = \frac{1}{|\mathcal{X}^P\cap\Phi^{k}|}\sum_{\mathbf{x}\in \mathcal{X}^P\cap\Phi^{k}}\bar{Y}(\mathbf{x})
	\end{equation} 
	Therefore, the objective function is the sum over all subregions the squared deviation from the mean. 
	Constraint~\eqref{eq:constrA} states that at least $N_{min}$ already sampled solutions must be clustered in one subregion of the partition, otherwise the optimization problem~\eqref{eq:min} is solved to minimum by a partition in which each subregion contains exactly one sample and the objective function value will be zero. An alternative is to limit the number of subregions to be divided into, denoted $d$ (Constraint~\eqref{eq:constrB}).
	The partition derived by solving problem~\eqref{eq:min} with constraints on the partition structure (e.g., $N_{min}$ or/and $d$), denoted $\mathcal{P}^*(\mathcal{R}^k)$, is one that optimally groups sampled points with similar performances. 
	Next, we discuss the solution algorithm for the proposed minimization problem~\eqref{eq:min}-\eqref{eq:constrB}.
	
	The proposed partition problem~\eqref{eq:min}-\eqref{eq:constrB} is similar to the underlying optimization problem of the decision tree model (e.g., \citet{breiman1984,dunn2018optimal}) with decision variables as input variables (features) and performances as dependent variables. 
	The decision tree model is mainly used for prediction purposes.
	The underlying optimization problem of training a decision tree is to find both a tree-structured split of the training samples based on input variables and labels for classification (or constant values for regression) of leaf nodes so that the prediction error is minimized, given some constraints on the tree structures (e.g., minimum leaf size and/or maximum tree depth) and restrictions against overfitting.
	The decision tree method CART by \citet{breiman1984} is a top-down greedy algorithm that does not guarantee a globally optimal decision tree. 
	A recent work of \citet{dunn2018optimal} proposes an algorithm that solves the decision tree problem at least locally optimally, and by initiating with multiple starting points, it attempts for global optimality. 
	They formulate the decision tree problem as an mixed integer programming problem and the objective function can be written as follows: 
	\begin{alignat}{3}
	\min_{\mathcal{T}} &\quad&  R(\mathcal{T})+\alpha|\mathcal{T}| & \label{eq:minOptimalTree}\\
	\text{s.t. } &\quad&  N(\ell)\geq N_{min}, & \quad \forall \quad\ell\in \text{leaves}(\mathcal{T}), \label{eq:constrOptimalTree}\\
	&\quad& depth(\mathcal{T}) \leq D, \label{eq:constrOptimalTreeB}
	\end{alignat}
	where $R(\mathcal{T})$ represents the prediction error of the tree $\mathcal{T}$ made on the training data. $|\mathcal{T}|$ denotes the number of branch node in tree $\mathcal{T}$. It is a representation of the tree complexity.  
	$N(\ell)$ represents the number of samples in each leave node $\ell$ of the tree $\mathcal{T}$, and hence the constraint~\eqref{eq:constrOptimalTree} restricts each leaf node to contain at least $N_{min}$ samples.
	One can also restricts the maximum depth of the tree (Constraint~\eqref{eq:constrOptimalTree}). It sets an upper bound on the number of leaf nodes (i.e., $2^{D}$). 
	In other words, the goal is to find a decision tree model that balances the prediction accuracy and the model complexity. The optimal decision tree method is shown to perform robustly to noisy training data (both feature and label noises).
	
	Note that a given tree $\mathcal{T}$ creates a unique partitioning on the feature space and hence a unique division of the samples. For regression tasks, the constant prediction value of a leaf node that minimizes the mean squared error is the mean performance of the samples in that leaf node.
	Thus, by setting $\alpha=0$ and choosing the mean squared error as the error metric (i.e., $R(\cdot)$), we retrieve an optimization problem that is equivalent to ~\eqref{eq:min}. 
	In other words, the proposed partitioning problem~\eqref{eq:min} is equivalent to the optimal parallel regression tree optimization \citep[Chapter 4]{dunn2018optimal} with complexity penalty coefficient $\alpha=0$. 
	The optimization problem~\eqref{eq:minOptimalTree}-\eqref{eq:constrOptimalTreeB} can be solved efficiently with their local search algorithm coupled with multiple start points for sample sizes up to hundreds of thousands. This is more than enough for discrete SO problems, in which samples are usually more time-consuming to simulate and hence limited in size. 
	
	\subsection{Hyperplane partition}\label{sec:hyperplanePartition}
	A hyperplane partition is one that contains cuts of the form $\mathbf{a}^T\mathbf{x}<b$ or $\mathbf{a}^T\mathbf{x}\geq b$ where $a_i \in \mathbb{R}$. The subregions resulting from a hyperplane partition are polyhedrons. The sampling strategy MIX-D algorithm \citep{pichitlamken2003combined} can efficiently sample uniformly from such subregion. 
	Hence, although more complicated than parallel partition, this family of partitions are also practical, especially when problem-specific structures of such form is known a priori.
	
	Assume that variables $y_j = \mathbf{a}_j^T\mathbf{x}$ for $j=1,...,q$ are known in advance to be potentially important splitting factors other than the decision variables $x_i$ for $i = 1,...,p$.
	Let us take the car-sharing fleet assignment problem as an example. The car-sharing service provider wants to find a fleet assignment to stations across the network that maximizes their profit. Geographically nearby stations often share demands among one another (i.e., spatial dependence), since customers are often willing to walk short distances for available vehicles if their target station runs out of vehicle. Thus, the total number of vehicles in the cluster of nearby stations is potentially a more important factor that influences the profit generated than the number of vehicles in each individual station. 
	Therefore, a partition based on the total number of vehicles in the clusters of nearby stations can more effectively divide the feasible region and lead to a more effective search for subregions with higher profits. This can result in a further improved finite-time performances and algorithm efficiency. 
	In this example, $y_j$ is the total number of vehicles assigned in the cluster of nearby stations $C_j$ as $y_j = \sum_{i\in C_j} x_i$ where $x_i$ are the number of vehicle assigned to station $i$. 
	
	The modification in order to incorporate such a hyperplane cut is minor. We simply treat $y_j$ as candidates that the partitioning algorithm can split on together with all the decision variables $x_i$. The proposed partitioning strategy by solving the optimization~\eqref{eq:min}-\eqref{eq:constrB} is robust to correlated splitting features such as $y_j$ and $x_i$, since it naturally includes a splitting feature selection process. 
	One other benefit the proposed adaptive partitioning strategy brings is that each subregion constructed contains at least one previously sampled solution, which can be used directly to initiate the MIX-D sampling scheme. 
	
	\subsection{Adaptive partitioning ESB\&B algorithm}
	The proposed adaptive partitioning ESB\&B algorithm is given by Algorithm \ref{alg:mESBB}. 
	The parameters related to partitioning strategy that need to be specified are the maximum tree depth $D$ and/or minimum number of points each subregion contains $N_{min}$, both are related to the number of subregions the current best subregion can be divided into. 
	The algorithm chooses the number of subregions to divide the region into optimally, given the predetermined parameters $D$ and/or $N_{min}$.
	Different from other generic partitioning strategies, the derived partition from the proposed adaptive partitioning strategy ensures that each new subregions contains at least one sampled points, or at least $N_{min}$ if this parameter is stated. This benefits the MIX-D sampling scheme in the sense that these interior sampled points can be used directly to initiate random walk. 
	The termination of the algorithm is usually when the simulation budget is exhausted. We select the final solution $\mathbf{x}^*$ as the one with the maximum cumulative sample average.
	
	\begin{algorithm}		
		\caption{Adaptive partitioning ESB\&B algorithm \label{alg:mESBB}}
		\small
		\begin{enumerate}
			\item[0.] Initialization:
			\begin{enumerate}
				\item set iteration counter $k = 0$, initial partition $\mathcal{P}_0 = \{\mathcal{X}\}$, best subregion $\mathcal{R}^0 = \mathcal{X}$
				\item sample uniformly at random the initial $n_0$ solutions in the feasible region $\mathcal{X}$, simulate $\Delta n_F$ replications of each sample, and record them in set $\Phi^0$.
				\item set training set $\Psi^0=\Phi^0$
			\end{enumerate} 		
			\item[1.] Partitioning:\\ If the best subregion $\mathcal{R}^k$ is singleton:
			\begin{enumerate}
				\item set the new full partition $\mathcal{P}^{'}_k = \mathcal{P}_k$			
			\end{enumerate}
			else:
			\begin{enumerate}
				\item construct a partition of the best subregion $\mathcal{P}(\mathcal{R}^k)$ using the proposed adaptive partitioning strategy with training set $\Psi^k$, splitting features $x_i$, predetermined parameters $D$ and/or $N_{min}$ and/or additional splitting features $y_j$
				\item define the new full partition by $\mathcal{P}^{'}_k = (\mathcal{P}_k\setminus \{\mathcal{R}^k\})\cup\mathcal{P}(\mathcal{R}^k)$
				\item denote $\mathcal{X}^P$ the elements of $\mathcal{P}^{'}_k$
			\end{enumerate}
			\item[2.] Sampling and bounding:
			\begin{enumerate}
				\item[2.1] Solution sampling: 
				\begin{enumerate}
					\item for each subregion $\mathcal{X}^P\in\mathcal{P}(\mathcal{R}^k)$, randomly sample $\nu_R$ solutions
					\item if $k > 0$ ,for each subregion $\mathcal{X}^P\in \mathcal{P}_k\setminus \{\mathcal{R}^k\}$, sample $\theta(\mathcal{X}^P)$ solutions, where $\theta(\mathcal{X}^P)$ is computed in Step 2.3 of iteration $k-1$ based on information in $\Phi^{k-1}$
					\item aggregate all of the sampled solutions into set $S^k$, and set $\Phi^k = \Phi^{k-1}\cup S^k$		
				\end{enumerate} 
				\item[2.2] Bound estimation:
				\begin{enumerate}
					\item for each $\mathbf{x}\in S^k$, simulate $\Delta n_F$ replications if $\mathbf{x}\notin \Phi^{k-1}$, simulate $\Delta n_A$ additional replications if $\mathbf{x}\in \Phi^{k-1}$
					\item for all $\mathcal{X}^P\in \mathcal{P}^{'}_k$, calculate estimates $\eta^{k+1}(\mathcal{X}^P)$
				\end{enumerate}
				\item[2.3] Sample allocation: compute the number of solutions to be sampled, $\theta(\mathcal{X}^P)$, for all $\mathcal{X}^P\in \mathcal{P}^{'}_k$ for iteration $k+1$ based on information in $\Phi^k$
			\end{enumerate}
			\item[3.] Updating partition and best subregion:
			\begin{enumerate}
				\item update the best subregion $\mathcal{R}^{k+1} = \arg\max_{\mathcal{X}^P\in\mathcal{P}^{'}_k}\{\eta^{k+1}(\mathcal{X}^P)\}$
				\item partition $\mathcal{P}_{k+1} = \mathcal{P}^{'}_k$
				\item training set $\Psi^{k+1} = \{\mathbf{x}:\mathbf{x}\in \mathcal{R}^{k+1}\cap\Phi^k\}$
				\item set $k=k+1$, go to Step 1.
			\end{enumerate}
		\end{enumerate}
	\end{algorithm}
	
	\section{Numerical examples}\label{sec:evaluation}
	In this section, we consider two sets of numerical examples to illustrate the proposed adaptive partitioning ESB\&B algorithm. The first example, Griewank function, has many local minima which makes it a challenging test function. 
	This low-dimensional example illustrates how the proposed method can escape tricky local minima by adaptively partitioning the feasible region with problem structure inferred from previously sampled points. 
	The second example is a real-world car-sharing case study in \citet{zhouOso2019} for which the optimal solutions are not known, it is used to explore the performance of the proposed algorithm on high-dimensional discrete SO problems. 
	In this example, we also demonstrate the use of hyperplane partition resulting from additional splitting variables $y_j$, which are derived from prior knowledge on the geographical locations of the stations.
	For both examples, we benchmark the proposed algorithm against the original ESB\&B algorithm with the generic partitioning strategy.
	
	\subsection{The Griewank function}
	We consider the Griewank function, which is commonly used as a test case for optimization algorithms \citep{salemi2019gaussian}. This function has a single global minimum and many local minima, which makes it a challenging test for optimization algorithms. 
	\subsubsection{Global minimum at the center of the feasible space}
\begin{figure}[htp]
		\minipage{1\textwidth}
		\centering
		\includegraphics[width=1\linewidth]{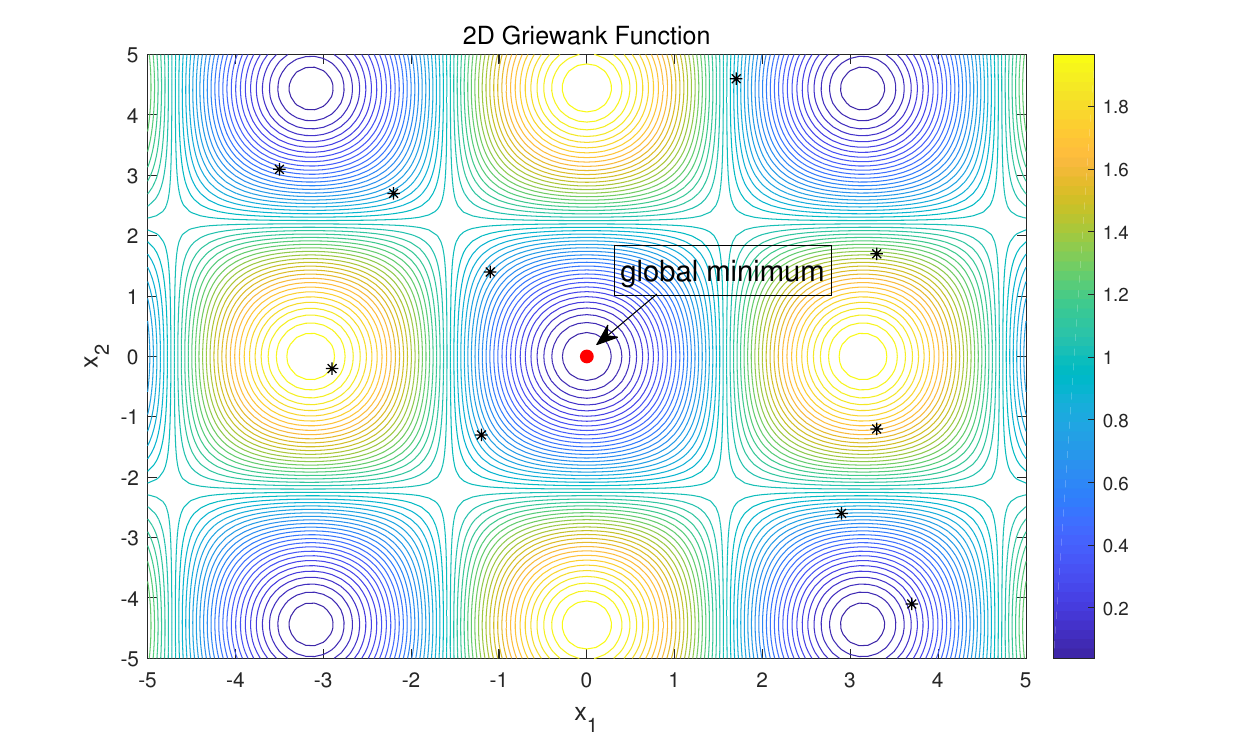}
		\caption{The contour plot of two-dimensional Griewank function on $[-5,5]\times[-5,5]$.}
		\label{fig:griewank_d2}
		\endminipage\hfill
	\end{figure}

	Figure \ref{fig:griewank_d2} displays the contour plot of the two-dimensional Griewank function on domain $[-5,5]\times[-5,5]$. The global minimum of the function is at the origin $(0,0)$ with response value $0$ and there are four local minima near the four corners of the domain with response values $0.0086$. 
	
	We first consider a minimization of the Griewank function with feasible region $[-5,5]\times[-5,5]$, where the globally optimal solution $(0,0)$ is at the center of the feasible region. To make it a discrete SO problem, the feasible region is divided into a $101\times101$ lattice, and the function value at each solution is given by the Griewank function plus a normally distributed noise with mean zero and variance $\sigma^2$. 
	In this numerical example, we take $\sigma = 0.01$, which is relative to the difference between local minima $0.0086$ and global minimum $0$. The number of replications for each non-encountered sample is set at $\Delta n_F=10$, and encountered sample $\Delta n_A = 2$.
	For both algorithms, we initiate with the same uniformly randomly sampled pool of size $10$ (asterisks in Figure \ref{fig:griewank_d2}). At each iteration, the total sampling budget for subregions other than the current best (denoted $\nu_O$) is $5$, and that for the current best subregion (denoted $|\mathcal{P}(\mathcal{R}^k)|\nu_R$) is $10$; the budget limit is $40$ iterations. Thus, at the end of each run, roughly $5\%$ of all feasible solutions are sampled and evaluated.
	For the generic partitioning scheme used in the original ESB\&B algorithm, the current best subregion $\mathcal{R}^k$ is divided into $\omega=2$ subregions equally along the longest dimension of $\mathcal{R}^k$. For the proposed adaptive partitioning scheme, the maximum tree depth $D$ is set at $2$, i.e., $\mathcal{R}^k$ can be at most divided into $4$ subregions, and the minimum number of sampled points grouped in one subregion is set at $2$, i.e., $N_{min} = 2$. We run each algorithm $50$ times.
	
\begin{figure}[htp]
		\minipage{1\textwidth}
		\centering
		\includegraphics[width=0.9\linewidth]{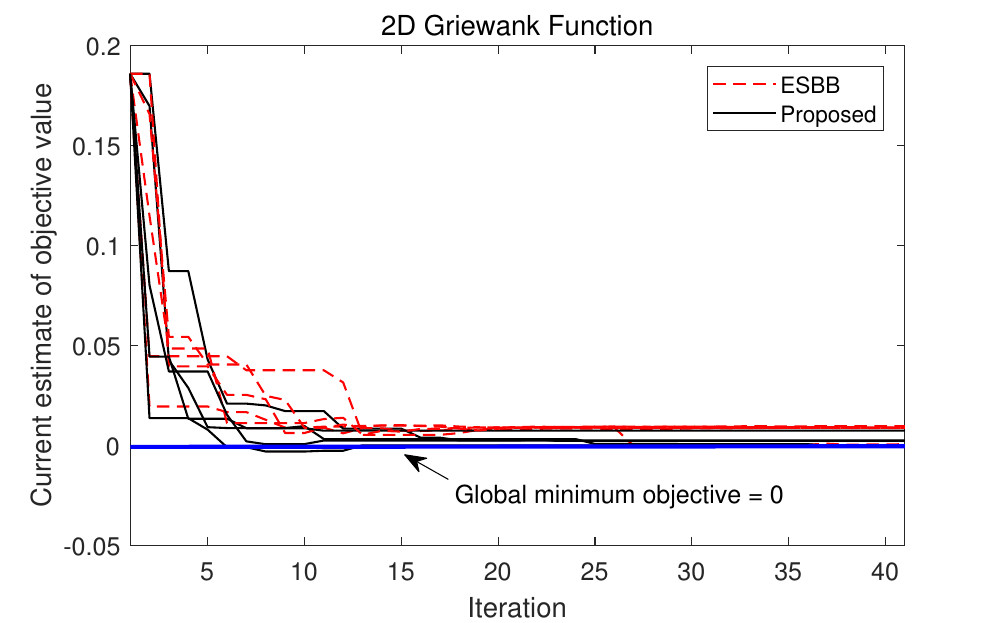}
		\caption{Objective function estimate of the current iterate across iterations.}
		\label{fig:GF_perform_sto}
		\endminipage\hfill
	\end{figure}

\begin{figure}[htp]
		\minipage{1\textwidth}
		\centering
		\includegraphics[width=0.9\linewidth]{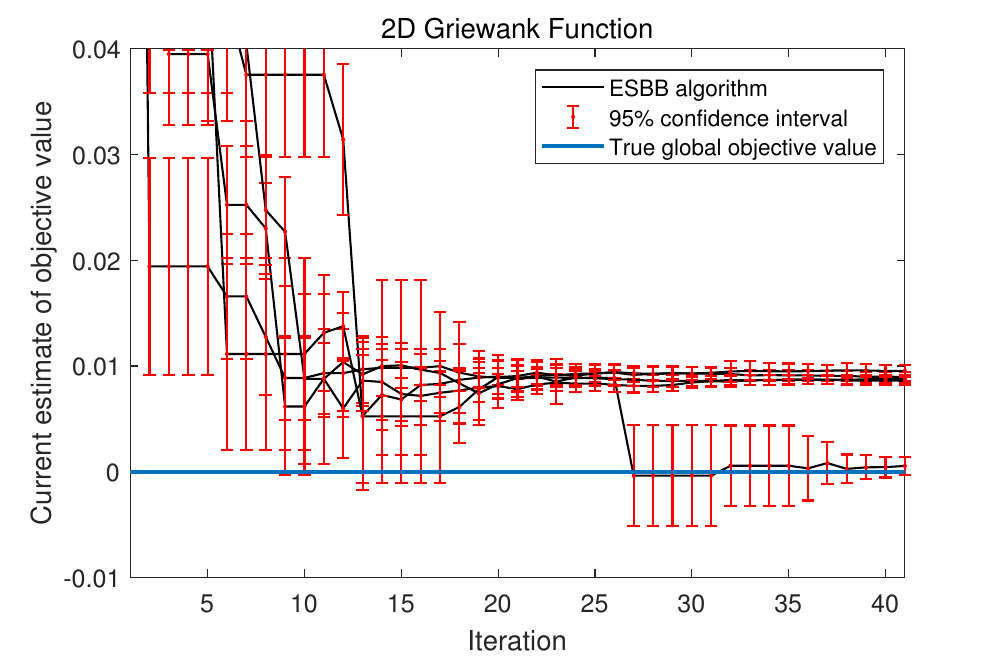}
		\caption{Objective function estimate of the current iterate with 95\% confidence interval across iterations of ESBB algorithm (zoomed-in results).}
		\label{fig:GF_perform_zoom_sto_esbb}
		\endminipage\hfill
	\end{figure}

\begin{figure}[htp]
		\minipage{1\textwidth}
		\centering
		\includegraphics[width=0.85\linewidth]{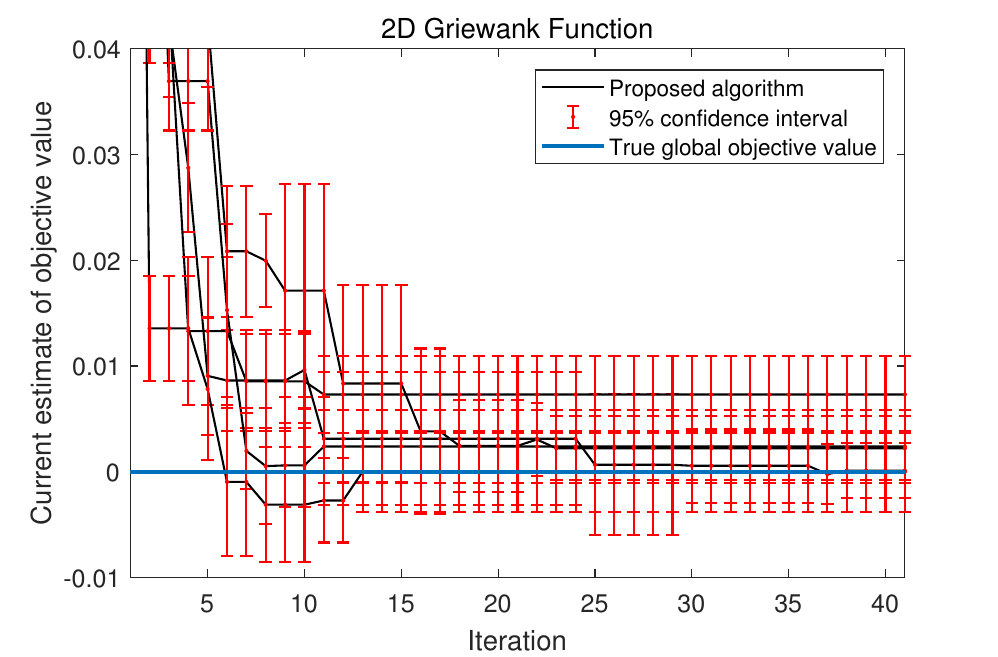}
		\caption{Objective function estimate of the current iterate with 95\% confidence interval across iterations of the proposed algorithm (zoomed-in results).}
		\label{fig:GF_perform_zoom_sto_otp}
		\endminipage\hfill
	\end{figure}

	Figure \ref{fig:GF_perform_sto} shows the current best estimate of the objective value across iterations of five randomly selected runs of each algorithm (i.e., five sample paths of each algorithm).
	The solid black lines display the results for the proposed algorithm, and the dashed red lines for the original ESB\&B algorithm. 
	As the iteration advances, the current best estimate of the objective value has a general decreasing trend for both algorithms, although there are temporary increases due to stochasticity. 
	Figure \ref{fig:GF_perform_zoom_sto_esbb} (resp. Figure \ref{fig:GF_perform_zoom_sto_otp}) plots a zoomed version of the current best estimate of the objective value with $95\%$ confidence intervals, so as to clearly show the performance of the ESB\&B algorithm (resp. proposed algorithm) close to the global minimum function value at zero (solid blue line). 
	Note that at each iteration, the current best solution is simulated at least $10$ replications, since each non-encountered sample will be simulated $\Delta n_F=10$ replications, and if it is sampled again, an additional $\Delta n_A = 2$ replications will be simulated.
	There is only one run of the ESB\&B algorithm that ends up with an estimated objective value that cannot be rejected at confidence level $95\%$ to be different from the true global minimum value $0$. 
	On the other hand, four runs of the proposed algorithm end up with estimated objective values that are statistically indifferent from the true global minimum after iteration $15$ at confidence level $95\%$.
	
\begin{figure}[htp]
		\minipage{1\textwidth}
		\centering
		\includegraphics[width=0.85\linewidth]{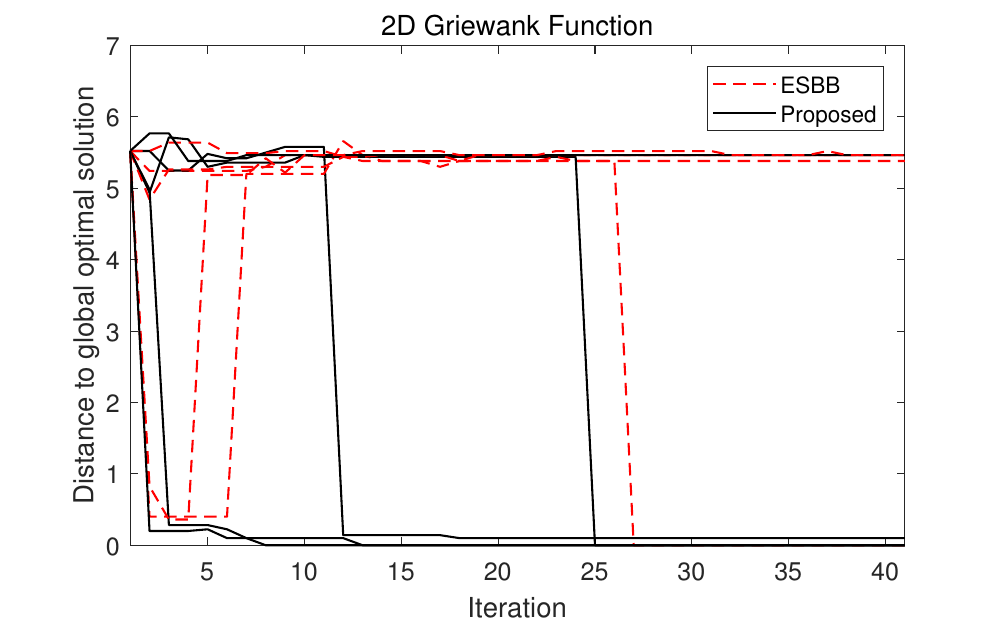}
		\caption{Distance between current best solution and the global minimum solution across iterations.}
		\label{fig:GF_dist_sto}
		\endminipage\hfill
	\end{figure}

	Figure \ref{fig:GF_dist_sto} shows the distance between current best solution and the global minimum solution across iterations. The solid black lines display the results for the proposed algorithm and dashed red lines for the ESB\&B algorithm. Four runs of the proposed algorithm end up with solutions close to the true global minimum at $(0,0)$, whereas only one run of the ESB\&B algorithm ends up close to $(0,0)$.
	Based on these five experiment runs, the proposed method tends to find solutions that are closer to the global minimum. 
	
\begin{figure}[htp]
		\minipage{1\textwidth}
		\centering
		\includegraphics[width=0.9\linewidth]{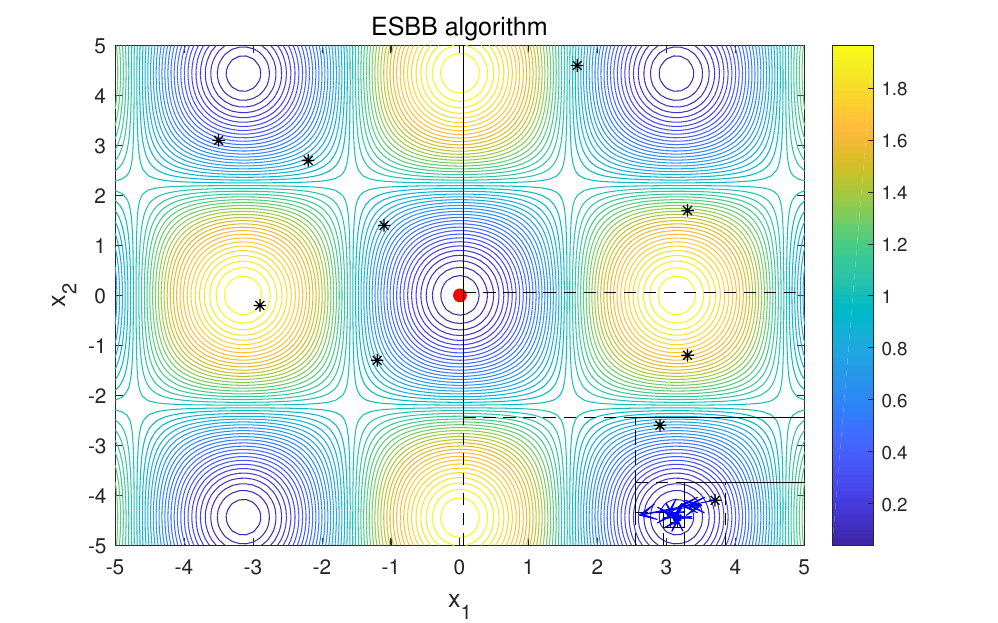}
		\caption{The path of best solution at current iterate across iterations in the feasible domain of the original ESB\&B algorithm.}
		\label{fig:ESBB_journey}
		\endminipage\hfill
		\minipage{1\textwidth}
		\centering
		\includegraphics[width=0.9\linewidth]{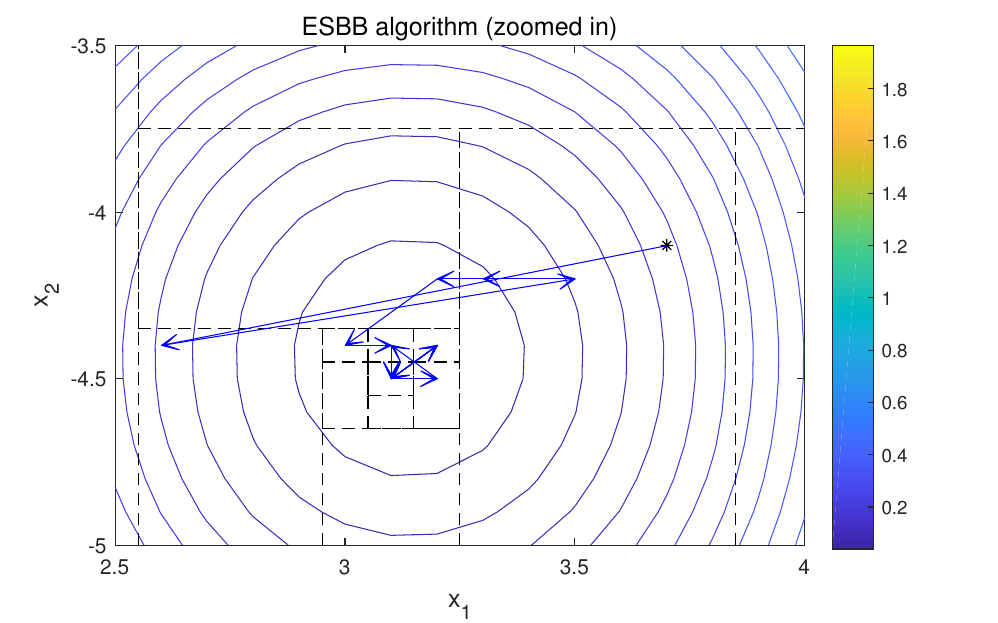}
		\caption{The path of best solution at current iterate across iterations in the feasible domain of the original ESB\&B algorithm (zoomed-in results).}
		\label{fig:ESBB_journey_zoom}
		\endminipage\hfill
	\end{figure}

\begin{figure}[htp]
		\minipage{1\textwidth}
		\centering
		\includegraphics[width=0.85\linewidth]{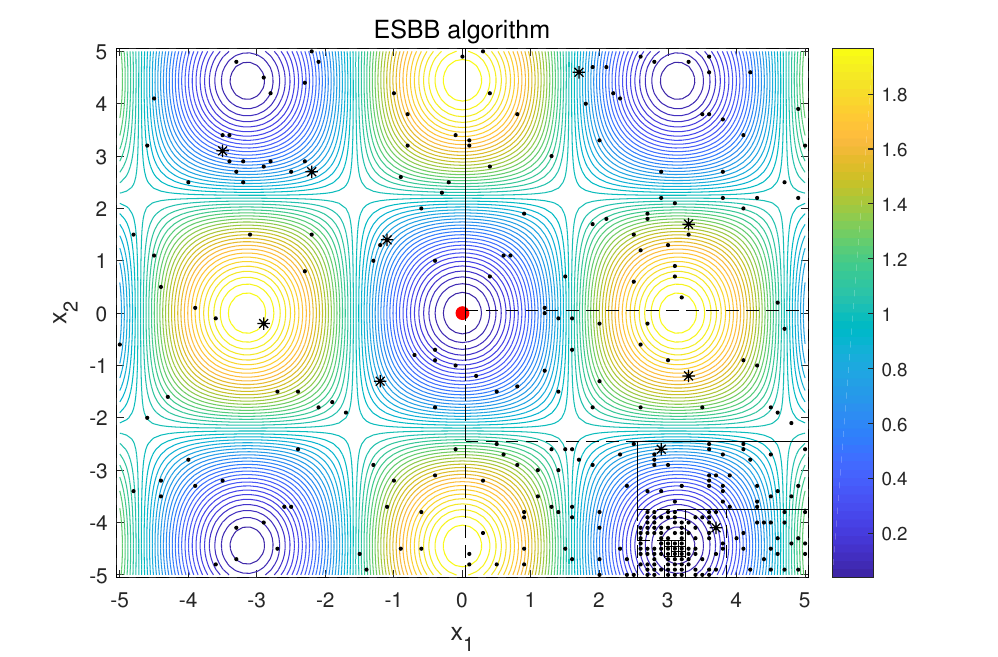}
		\caption{Allocation of sampling budget in the feasible domain of the original ESB\&B algorithm.}
		\label{fig:ESBB_sampling}
		\endminipage\hfill
	\end{figure}

	Figures \ref{fig:ESBB_journey} and \ref{fig:ESBB_sampling} show the experiment results for the first run of the original ESB\&B algorithm, which is a typical run of ESB\&B getting trapped at a local minimum solution.
	In both figures, the final partition of the feasible region and the contour plot of the Griewank function are displayed. 
	Figure \ref{fig:ESBB_journey} plots the path of the best solution at current iterate across iterations in the feasible domain, the path is plotted with blue arrows, and Figure \ref{fig:ESBB_journey_zoom} displays the zoomed-in results. The asterisk points plotted are the initial sampled set $\Phi^0$.
	The generic partitioning strategy used in the ESB\&B algorithm missed the subregion that contains the global minimum in the first iteration and left it near the boundary of a subregion. The algorithm gets trapped to the lower right local minimum. 
	Figure \ref{fig:ESBB_sampling} displays the sampling budget allocation in the feasible domain at the end of iteration $40$, in which each black dot represents a sampled solution. As discussed, much less sampling budget is allocated to subregions that have been discarded earlier.
	
	\begin{figure}[htp]
		\minipage{1\textwidth}
		\centering
		\includegraphics[width=0.85\linewidth]{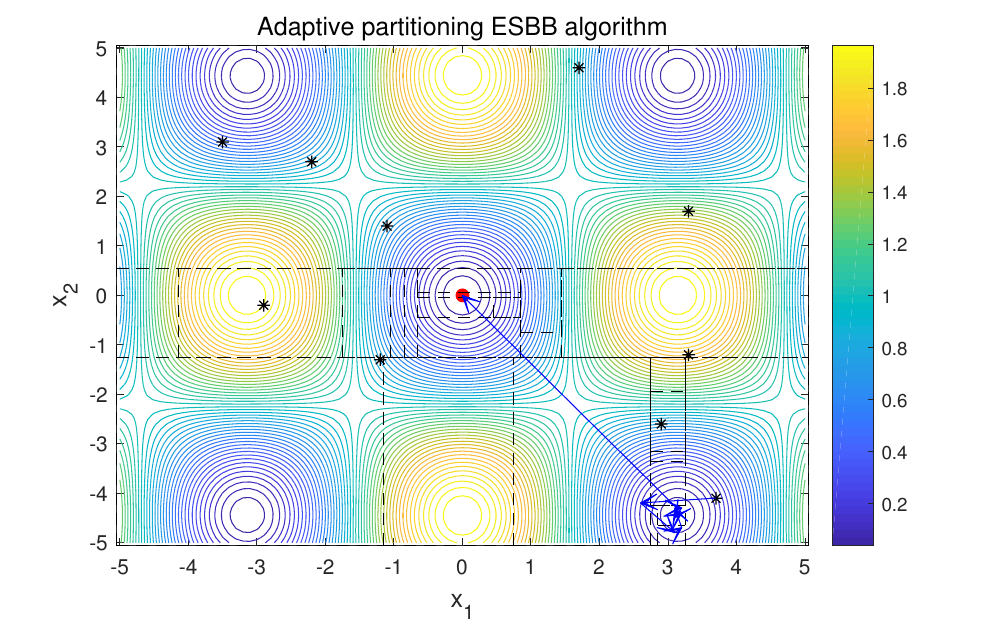}
		\caption{The path of best solution at current iterate across iterations in the feasible domain of the proposed algorithm.}
		\label{fig:OPT_journey}
		\endminipage\hfill
	\end{figure}
	
\begin{figure}[htp]
		\minipage{1\textwidth}
		\centering
		\includegraphics[width=0.9\linewidth]{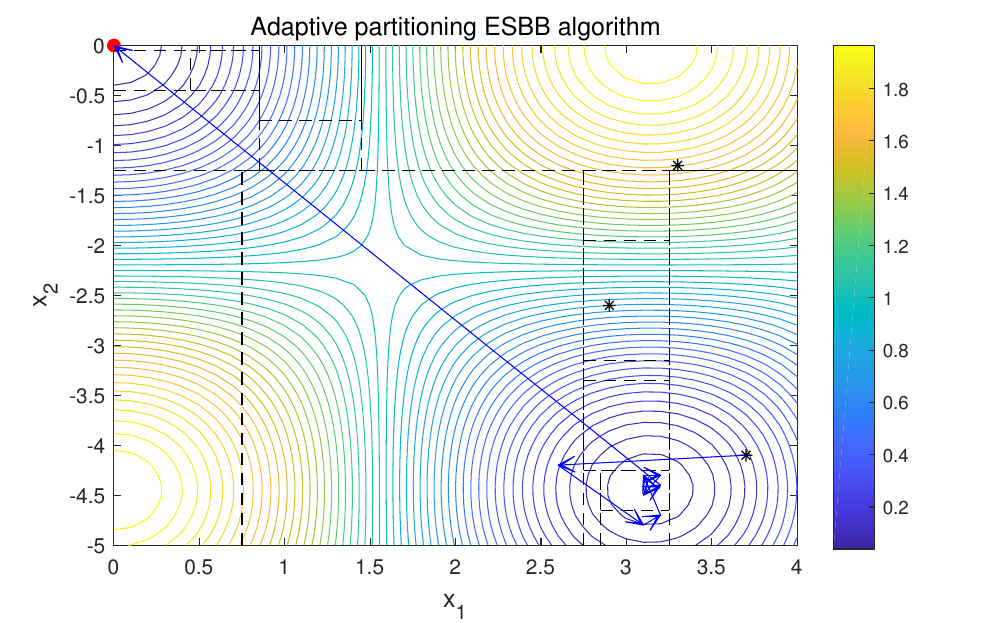}
		\caption{The path of best solution at current iterate across iterations in the feasible domain of the proposed algorithm (zoomed-in results).}
		\label{fig:OPT_journey_zoom}
		\endminipage\hfill
	\end{figure}

\begin{figure}[htp]
		\minipage{1\textwidth}
		\centering
		\includegraphics[width=0.9\linewidth]{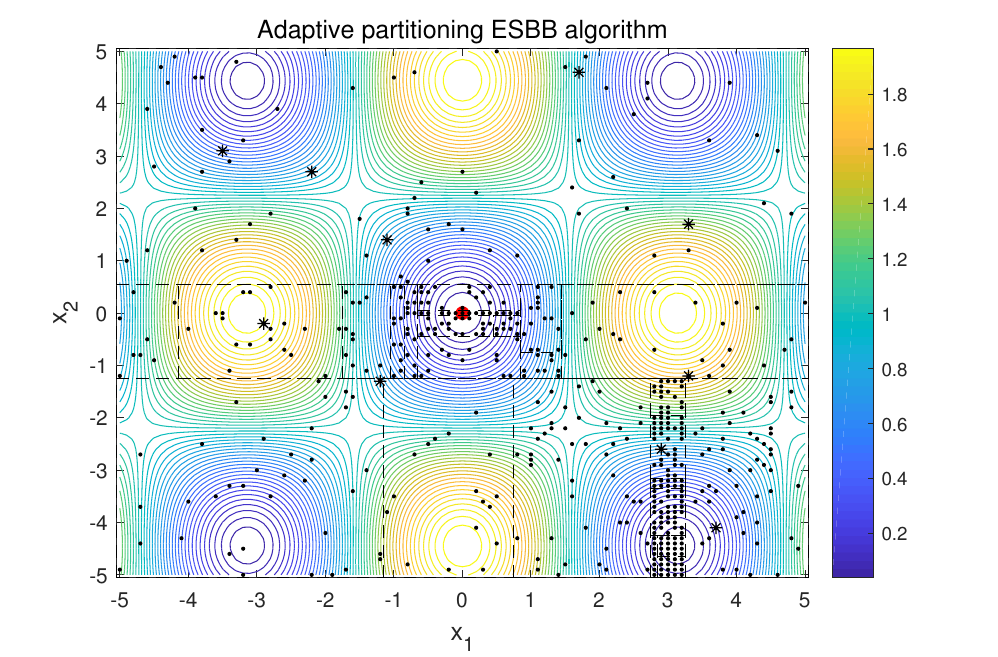}
		\caption{Allocation of sampling budget in the feasible domain of the proposed algorithm.}
		\label{fig:OPT_sampling}
		\endminipage\hfill
	\end{figure}

	Figures \ref{fig:OPT_journey} and \ref{fig:OPT_sampling} show the experiment results for first run of the proposed algorithm. 
	As before, the final partition of the feasible region and the contour plot of the Griewank function are displayed in both figures. 
	Figure \ref{fig:OPT_journey} plots the path of the current best solution across iterations in the feasible domain, and Figure \ref{fig:OPT_journey_zoom} displays the zoomed-in results.
	Figure \ref{fig:OPT_sampling} displays the sampling budget allocation in the feasible domain. These figures have similar layouts as Figure \ref{fig:ESBB_journey} and \ref{fig:ESBB_sampling}, respectively.
	Different from the generic partitioning strategy, the proposed adaptive partitioning strategy divides the feasible region along $x_2$ initially and note that the globally optimal solution is placed in the interior of one of the subregion. At iteration 24, the algorithm escapes the lower right local minimum and starts to explore the middle subregion. The final partition of the feasible region identifies the underlying function's peaks and basins successfully. The sampling budgets are allocated more to the subregions that contain the basins than those containing the peaks as expected.

	\begin{figure}[htp]
		\minipage{1\textwidth}
		\centering
		\includegraphics[width=0.9\linewidth]{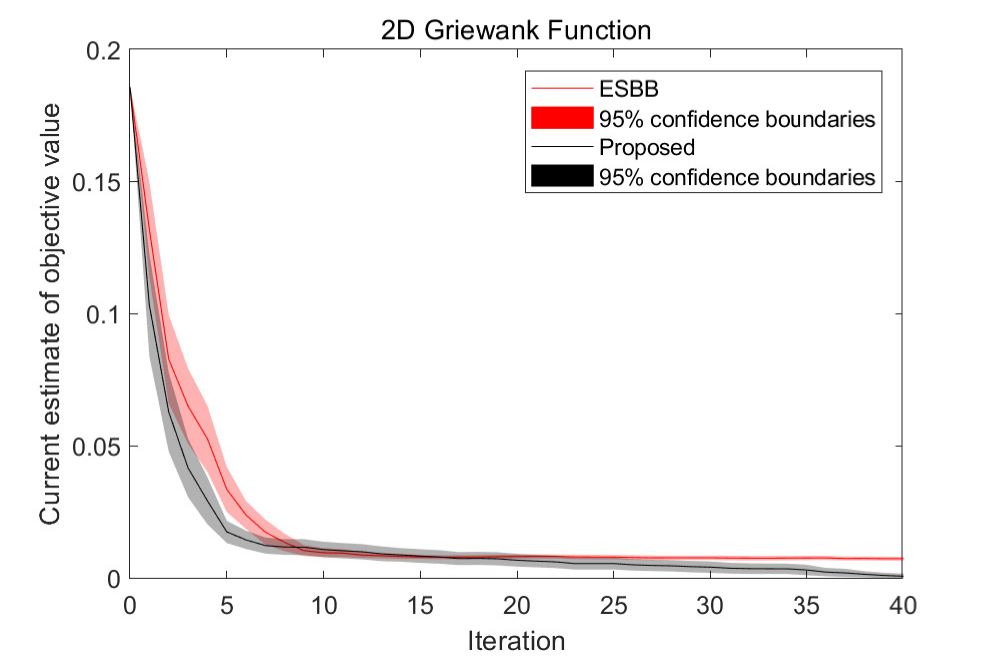}
		\caption{Objective function estimate of the current iterate across iterations averaged over 50 algorithm runs.}
		\label{fig:GF_perform_sto_agg}
		\endminipage\hfill
		\minipage{1\textwidth}
		\centering
		\includegraphics[width=0.9\linewidth]{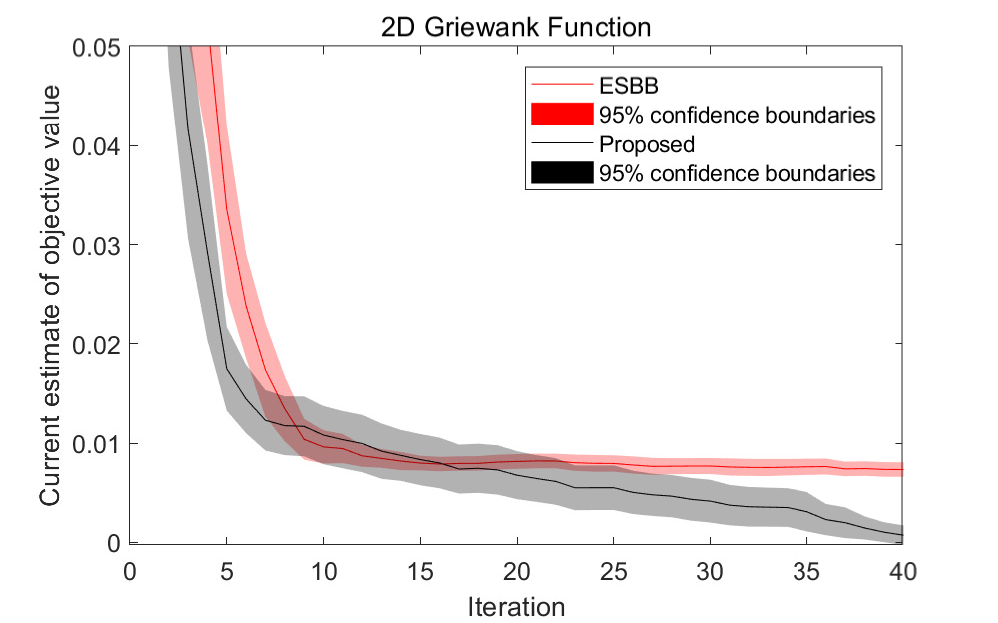}
		\caption{Objective function estimate of the current iterate across iterations averaged over 50 algorithm runs (zoomed-in results).}
		\label{fig:GF_perform_zoom_sto_agg}
		\endminipage\hfill
	\end{figure}
	
	An average sample-path performance (over the 50 runs) for each algorithm is constructed; this performance metric is used to compare against other algorithms, e.g., \citet{xu2010industrial,xu2013adaptive,xu2013empirical}.
	Figure \ref{fig:GF_perform_sto_agg} plots the objective value of the estimated optimal solution across iterations (averaged over 50 algorithm runs) for the ESB\&B algorithm (red line with shaded $95\%$ confidence boundary) and the proposed algorithm (black line with shaded $95\%$ confidence boundary). 
	Figure \ref{fig:GF_perform_zoom_sto_agg} zooms in on the average performance close to the global minimum function value at zero. 
	Since this is a minimization problem, the lower the curve, the better the algorithm performance in terms of estimated objective value. Note that the curve for the proposed algorithm is lower than that for the ESB\&B algorithm except for iterations $10$ to $25$ where they overlap with each other. This indicates: (i) at early stage, the proposed algorithm find solutions with improved performances faster; (ii) during iterations $10$ to $25$, both algorithm reaches solutions with similar performances around local optimal objective values; (iii) as the iteration number increases, the proposed algorithm continues to find solutions with improved performances whereas the ESB\&B algorithm mostly gets stuck around the locally optimal solutions. 
	At the termination of the algorithms (i.e., iteration $40$), the proposed algorithm ends up with an average estimated objective value that is statistically lower than that of the ESB\&B algorithm. 
	Additionally, the proposed algorithm finds the true globally optimal solution (i.e., $(0,0)$) for $27$ out of $50$ runs, whereas the ESB\&B algorithm finds it for $3$ out of $50$ runs. 
	The propose algorithm finds a final solution with mean performances statistically indifferent from the true optimal value (i.e., $0$) for $41$ out of $50$ runs at significance level $0.05$, whereas the ESB\&B algorithm finds it for $4$ out of $50$ runs.

	\subsubsection{Global minimum away from center of the feasible space}
\begin{figure}[htp]
		\minipage{1\textwidth}
		\centering
		\includegraphics[width=1\linewidth]{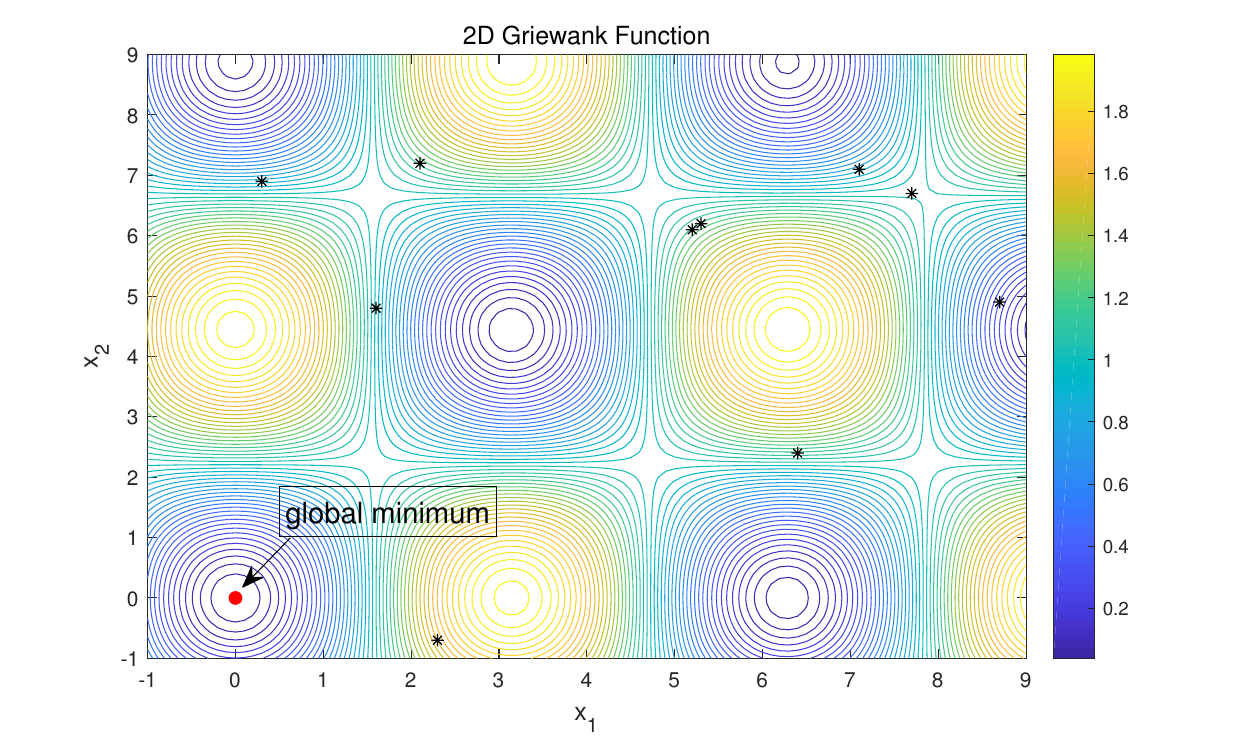}
		\caption{The contour plot of two-dimensional Griewank function on $[-1,9]\times[-1,9]$.}
		\label{fig:griewank_d2_shift}
		\endminipage\hfill
	\end{figure}

	Next, we consider a minimization of the Griewank function with feasible region $[-1,9]\times[-1,9]$, where the globally optimal solution $(0,0)$ is no longer at the center of the feasible region. In this example, the generic partitioning strategy of ESB\&B algorithm does not end up with an initial cut close to globally optimal solution. 
	Figure \ref{fig:griewank_d2_shift} displays the contour plot of the two-dimensional Griewank function on domain $[-1,9]\times[-1,9]$ together with a uniformly random generated initial sample set (displayed in asterisks) for both algorithms. Parameters of the algorithms are set the same as in the previous experiment. We run each algorithm $50$ times.
	
\begin{figure}[htp]
		\minipage{1\textwidth}
		\centering
		\includegraphics[width=0.9\linewidth]{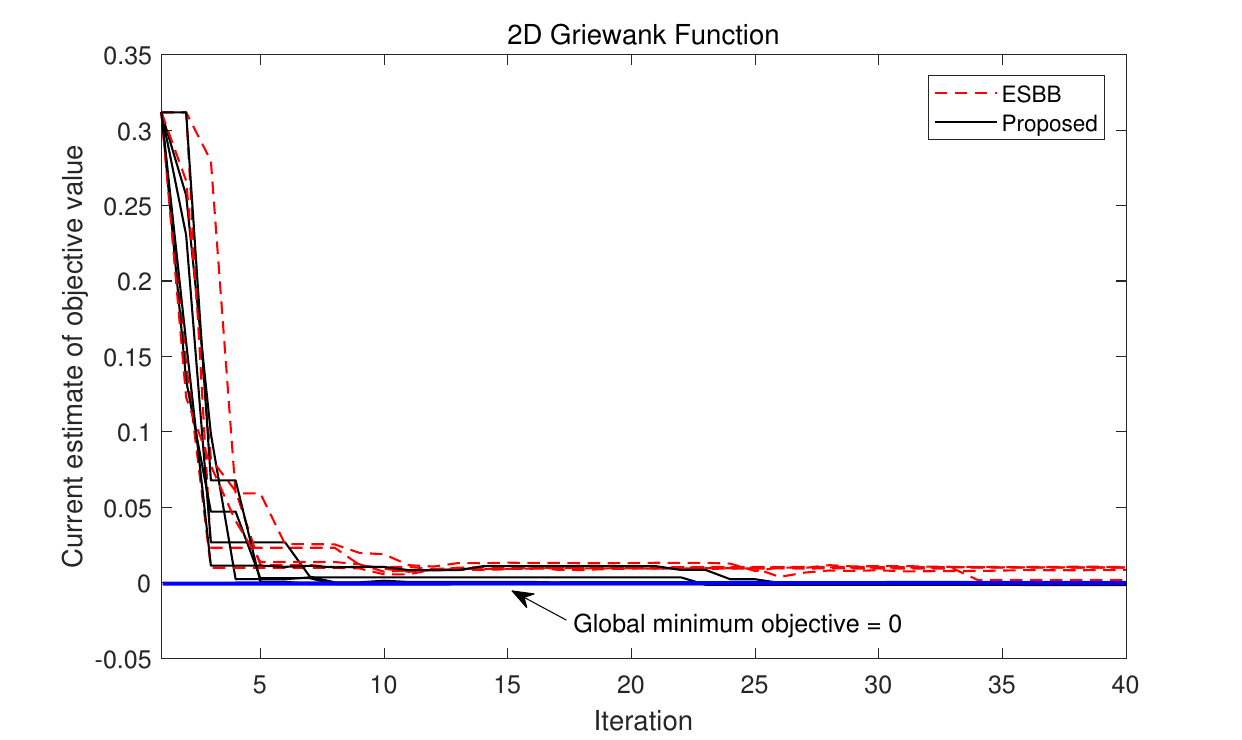}
		\caption{Objective function estimate of the current iterate across iterations.}
		\label{fig:GF_perform_sto_shift}
		\endminipage\hfill
	\end{figure}

\begin{figure}[htp]
		\minipage{1\textwidth}
		\centering
		\includegraphics[width=0.9\linewidth]{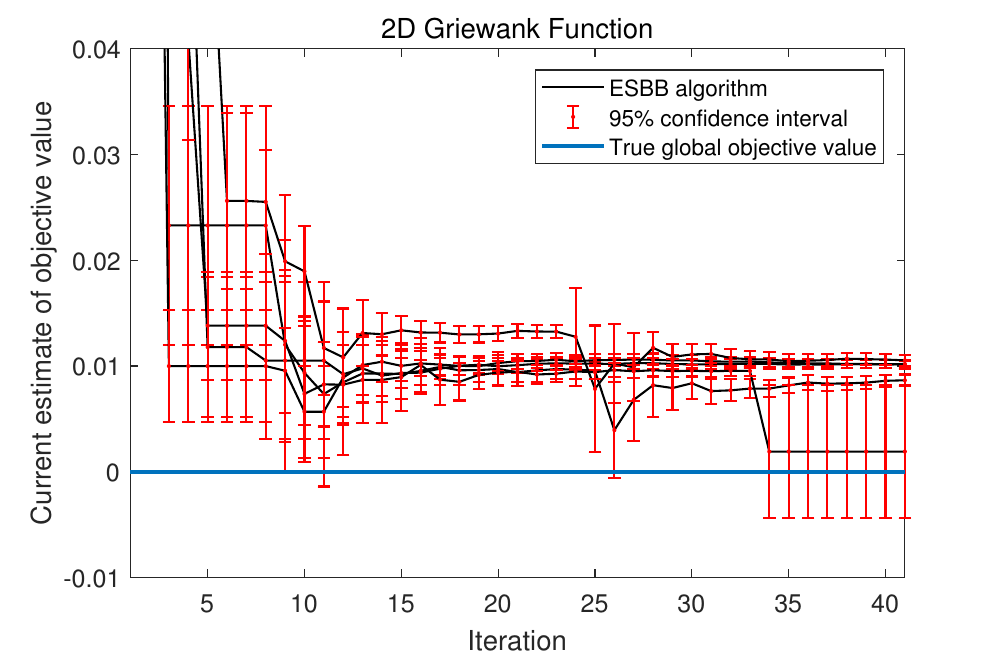}
		\caption{Objective function estimate of the current iterate with 95\% confidence interval across iterations of ESBB algorithm (zoomed-in results).}
		\label{fig:GF_perform_zoom_sto_shift_esbb}
		\endminipage\hfill
	\end{figure}

\begin{figure}[htp]
		\minipage{1\textwidth}
		\centering
		\includegraphics[width=0.85\linewidth]{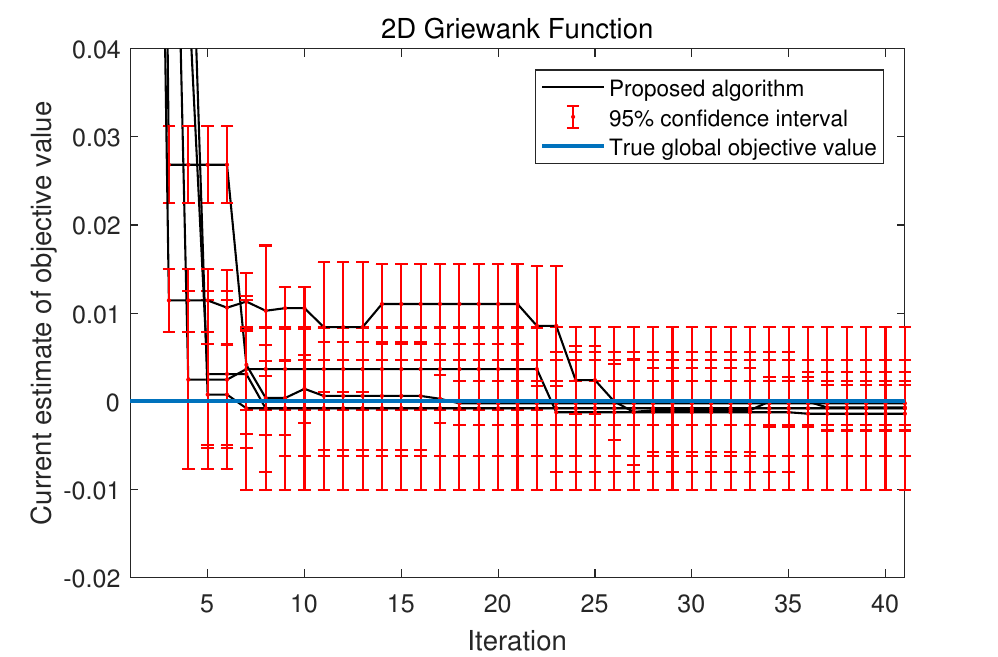}
		\caption{Objective function estimate of the current iterate with 95\% confidence interval across iterations of the proposed algorithm (zoomed-in results).}
		\label{fig:GF_perform_zoom_sto_shift_otp}
		\endminipage\hfill
	\end{figure}

	Figure \ref{fig:GF_perform_sto_shift} plots the current best estimate of objective value across iterations of five randomly selected runs of each algorithm (i.e., five sample paths of each algorithm). As before, we observe that the current best estimate of the objective value inherits a decreasing trend for both algorithms.
	Figure \ref{fig:GF_perform_zoom_sto_shift_esbb} (resp. Figure \ref{fig:GF_perform_zoom_sto_shift_otp}) plots the current best estimate of objective value with $95\%$ confidence interval of ESB\&B algorithm (resp. proposed algorithm) and zooms in on the performance close to the global minimum function value at zero (solid blue line). 
	There is only one run of the ESB\&B algorithm that ends up with estimated objective values that cannot be rejected at confidence level $95\%$ to be different from the true global minimum value $0$. 
	Meanwhile, all five runs of the proposed algorithm end up with estimated objective values that are statistically indifferent from the true global minimum after iteration $25$ at confidence level $95\%$.
	
\begin{figure}
		\minipage{1\textwidth}
		\centering
		\includegraphics[width=0.85\linewidth]{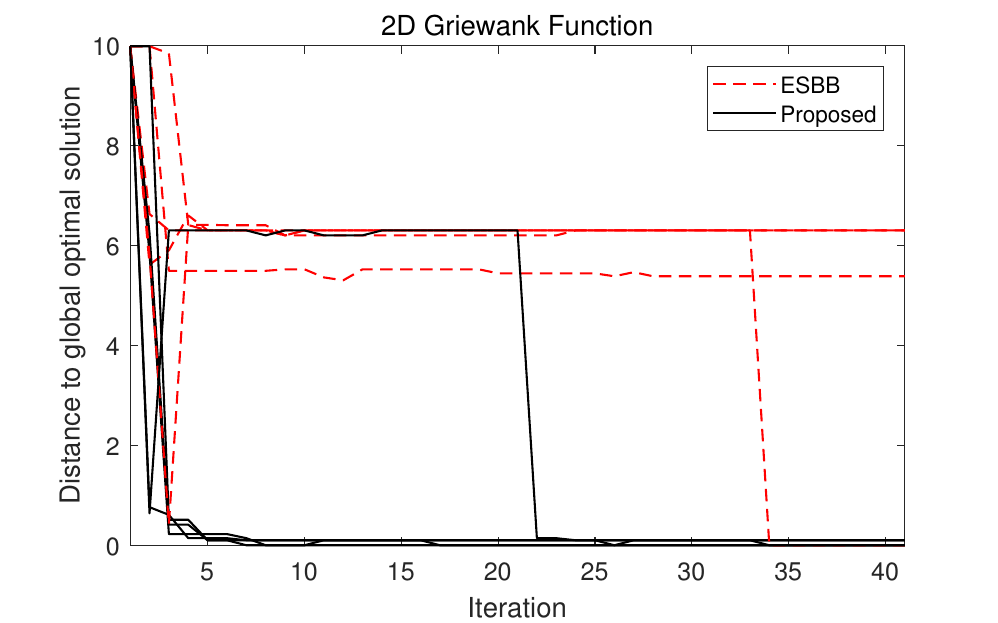}
		\caption{Distance between current best solution and the global minimum solution across iterations.}
		\label{fig:GF_dist_sto_shift}
		\endminipage\hfill
	\end{figure}

	Figure \ref{fig:GF_dist_sto_shift} shows the distance between current best solution and the global minimum solution across iterations. This figure has a similar layout as Figure \ref{fig:GF_dist_sto}. All runs of the proposed algorithm end up with solutions close to the true global minimum at $(0,0)$, whereas only one run of the ESB\&B algorithm ends up close to $(0,0)$. Based on these five experiment runs, the proposed methods tends to find solutions that are closer to the global minimum.

\begin{figure}[htp]
		\minipage{1\textwidth}
		\centering
		\includegraphics[width=0.9\linewidth]{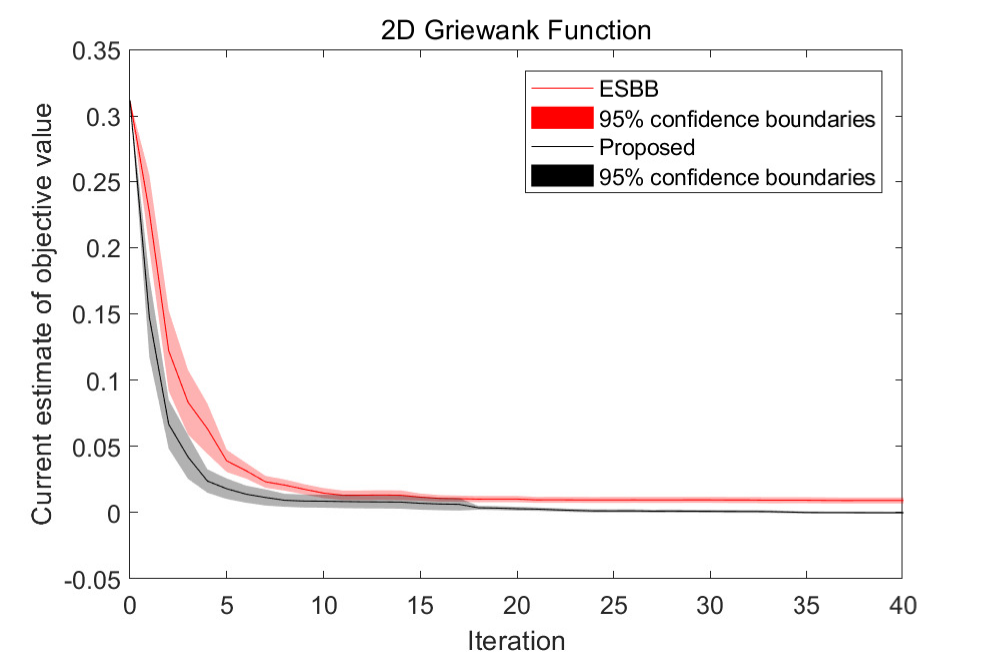}
		\caption{Objective function estimate of the current iterate across iterations averaged over 50 algorithm runs.}
		\label{fig:GF_perform_sto_shift_agg}
		\endminipage\hfill
		\minipage{1\textwidth}
		\centering
		\includegraphics[width=0.9\linewidth]{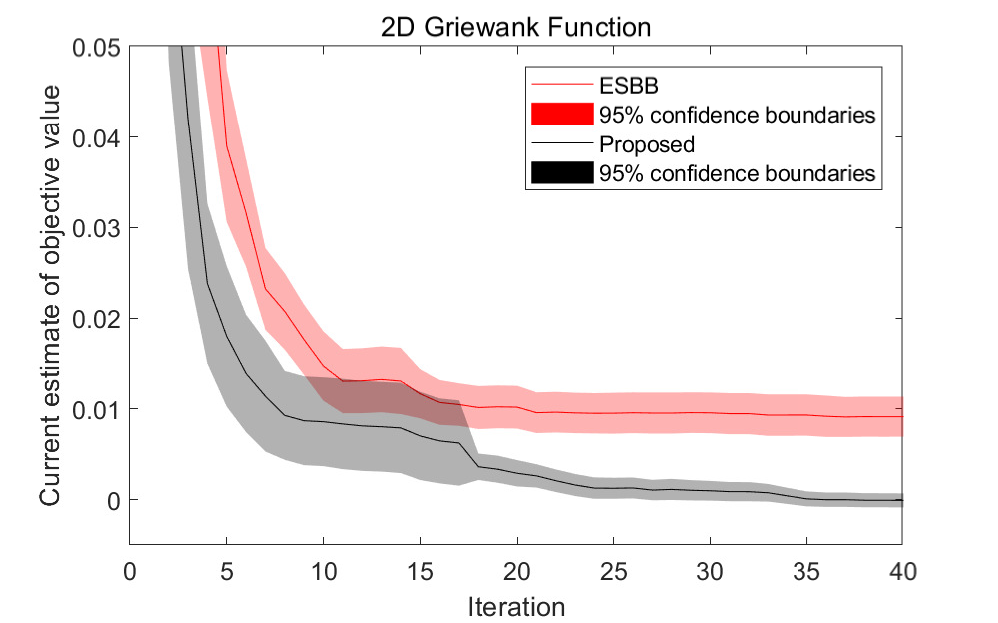}
		\caption{Objective function estimate of the current iterate across iterations averaged over 50 algorithm runs (zoomed-in results).}
		\label{fig:GF_perform_zoom_sto_shift_agg}
		\endminipage\hfill
	\end{figure}

	Figure \ref{fig:GF_perform_sto_shift_agg} plots the average estimate of the objective value across iterations for the ESB\&B algorithm (red line with shaded $95\%$ confidence boundary) and the proposed algorithm (black line with shaded $95\%$ confidence boundary). Figure \ref{fig:GF_perform_zoom_sto_shift_agg} zooms in on the average performance close to the global minimum function value of zero.  
	As before, the curve for the proposed algorithm is lower than that of the ESB\&B algorithm. In other words, the proposed algorithm has a better performance than the ESB\&B algorithm on average. 
	The initial decrease in estimated objective value of the proposed algorithm is faster than that of the ESB\&B algorithm, this indicates a faster exploring speed of solutions with improved performances. At termination (i.e., iteration 40), the proposed algorithm ends up with an average estimated objective value that is statistically lower than that of the ESB\&B algorithm. 
	Additionally, $35$ runs of the proposed algorithm end up with true globally optimal solution at $(0,0)$, whereas only $14$ runs for the ESB\&B algorithm obtain the true globally optimal solution. Moreover,
	$46$ runs of the proposed algorithm end with a solution with mean performance statistically indifference from the true optimal value $0$ at significance level $0.05$, whereas only $14$ runs for the ESB\&B algorithm achieve this.
	
	In summary, the proposed adaptive partitioning ESB\&B algorithm is a globally convergent algorithm, meanwhile it has an improved finite-time (limited sample budget) performance in terms of finding improved solutions faster, better estimated final objective values and higher likelihood of finding the true global optimal solutions. 
	
	\subsection{The car-sharing fleet assignment problem}
	In this section, we apply the proposed algorithm to a real-world car-sharing fleet assignment problem. This case study problem is adopted from \citet{zhouOso2019}.
	It considers a two-way car-sharing system from the perspective of the service operator. Essentially, it consists of finding an assignment of a fleet of vehicles across the network of stations that maximizes the expected profit over a given finite time horizon, denoted as the planning period. 
	Instead of using a simplified description of demand and of demand-supply interactions, this case study relies on a demand simulator of \citet{fields2017data} developed based on the rich high-resolution reservation data from Zipcar's Boston market. 
	This car-sharing fleet assignment problem is formulated as a discrete SO problem in \citet{zhouOso2019}, where the objective function value (i.e., the expected profit of a given fleet assignment) can only be obtained via the demand simulation.
	
	The discrete SO problem is formulated as follows:
	\begin{alignat}{2}
	\max_{\mathbf{x}} &\quad&&  E[R(\mathbf{x};\mathbf{q}_1)] - \sum_{i\in\mathcal{I}}c_ix_i  \label{eq:obj_zipcar}\\
	\text{s.t. } &\quad&&  \sum_{i\in\mathcal{I}} x_i\leq N  \label{eq:constr1_zipcar}\\
	&\quad&& x_i \leq N^i,  \quad \forall\quad i\in\mathcal{I} \label{eq:constr2_zipcar}\\
	&\quad&& x_i \in\mathbb{Z}_{\geq 0}, \quad \forall\quad i\in\mathcal{I} \label{eq:constr3_zipcar}
	\end{alignat}
	where $\mathbf{x} = [x_i]$ are the decision variables and $x_i$ is the number of vehicles assigned to station $i$. $R(\mathbf{x};\mathbf{q}_1)$ is the random variable representing the revenue with fleet assignment $\mathbf{x}$ and exogenous simulation parameter vector (e.g., reservation pricing) $\mathbf{q}_1$, $c_i$ is the exogenous cost, over the planning period, of a parking space at station $i$, $N$ is the total number of vehicles to assign, $N^i$ is the capacity of station $i$, and $\mathcal{I}$ is the set of all stations. 
	The objective function~\eqref{eq:obj_zipcar} represents the expected profit for a given fleet assignment $\mathbf{x}$ as the difference between the expected revenue $E[R(\mathbf{x};\mathbf{q}_1)]$ and costs $\sum_{i\in\mathcal{I}}c_ix_i$. The estimates of the expected revenue $E[R(\mathbf{x};\mathbf{q}_1)]$ can only be obtained via simulations. Constraint~\eqref{eq:constr1_zipcar} bounds the total number of vehicles assigned across all the stations. Constraint~\eqref{eq:constr2_zipcar} bounds the number of vehicles assigned to each individual station $i$ by the capacity of the station. Finally, the number of vehicles assigned to each station must be an non-negative integer (Constraint~\eqref{eq:constr3_zipcar}).
	
	As discussed in Section \ref{sec:hyperplanePartition}, in the car-sharing fleet assignment problem, the total number of vehicles assigned to a cluster of nearby stations (i.e., $y_j = \sum_{i\in C_j} x_i$) may form a more efficient cut of the feasible region, since stations nearby usually share demands as customers are likely to search nearby stations for substitutions if the target station runs out of vehicles. 
	Therefore, to address the formulated discrete SO problem~\eqref{eq:obj_zipcar}, we apply the following three algorithms: the original ESB\&B algorithm, the proposed algorithm with parallel partition, and the proposed algorithm with hyperplane partition.
	The potential splitting factors $y_j$ for hyperplane partition are generated as follows: for each stations $i$, form a cluster of stations center at station $i$ with a given radius, and eliminate replications. In this section, we choose the radius of $1$ in the distance unit used in the simulator when calculating the spillover effect from one station to another.
	
\begin{figure}[htp]
		\minipage{1\textwidth}
		\centering
		\includegraphics[width=0.7\linewidth]{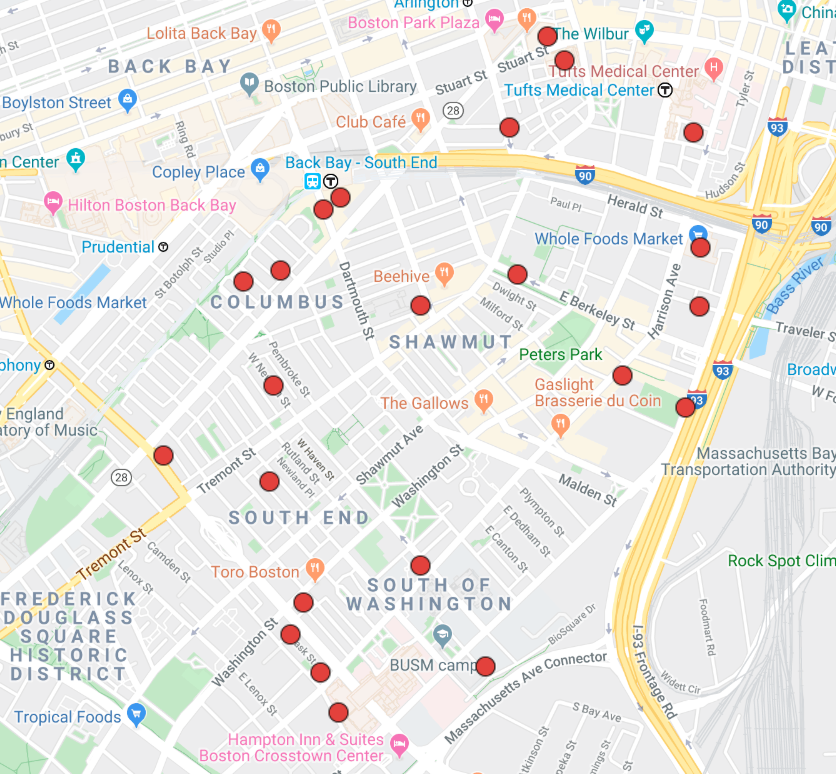}
		\caption{Zipcar stations in Boston South End neighborhood \citep{googlemapBos}}
		\label{fig:zipcar23}
		\endminipage\hfill
	\end{figure}

	We consider the fleet assignment problem of Boston south end area which contains $23$ stations, i.e., the decision vector is of dimension $23$ (shown in Figure \ref{fig:zipcar23}). Each station $i$ has a space capacity $N^i = 16$ and the total number of cars to assign is $N=211$. 
	All other exogenous variables (e.g., $c_i$) other than demand level of the simulator are set the same as in \citet{zhouOso2019}. We tested the algorithms on one low-demand level and one high-demand level case.
	The maximum number of algorithm iterations is set to $40$. At every iteration, the number of solutions to be sampled from subregions other than the current best subregion is set to $10$ (i.e., $\nu_O = 10$), and the total number of solutions to be sampled from the current best subregion is set to $20$ (i.e., $|\mathcal{P}(\mathcal{R}^k)|\nu_R = 20$). 
	The number of replications for each non-encountered sample is set at $\Delta n_F=5$, encountered sample $\Delta n_A = 2$. 
	For all algorithms, we initialize with $20$ randomly uniformly sampled solutions plus one solution with relatively good performance, which can be considered as a warm start. This warm-start solution is obtained by solving the analytical car-sharing fleet assignment model formulated as an Mixed Integer Programming (MIP) in \citet[Eq.(8)-(13)]{zhouOso2019}.
	
	For the original ESB\&B algorithm, the best subregion is divided into $3$ new subregions at each iteration (i.e., $\omega=3$). 
	For the proposed algorithm with both parallel partition and hyperplane partition, the maximum tree depth is set at $2$ (i.e., $D=2$), and hence the current best subregion can be divided into at most $4$ subregions. The minimum number of sampled points grouped in one subregion is set at $2$, i.e., $N_{min} = 2$. 
	Note that the computational time spend on simulation per iteration in this experiment setting is roughly $500$ seconds, the proposed partitioning strategy with both parallel and hyperplane partitions finishes within $2$ seconds, which is comparably trivial.
	We run each algorithm five times.
	
\begin{figure}[htp]
		\minipage{1\textwidth}
		\centering
		\includegraphics[width=0.9\linewidth]{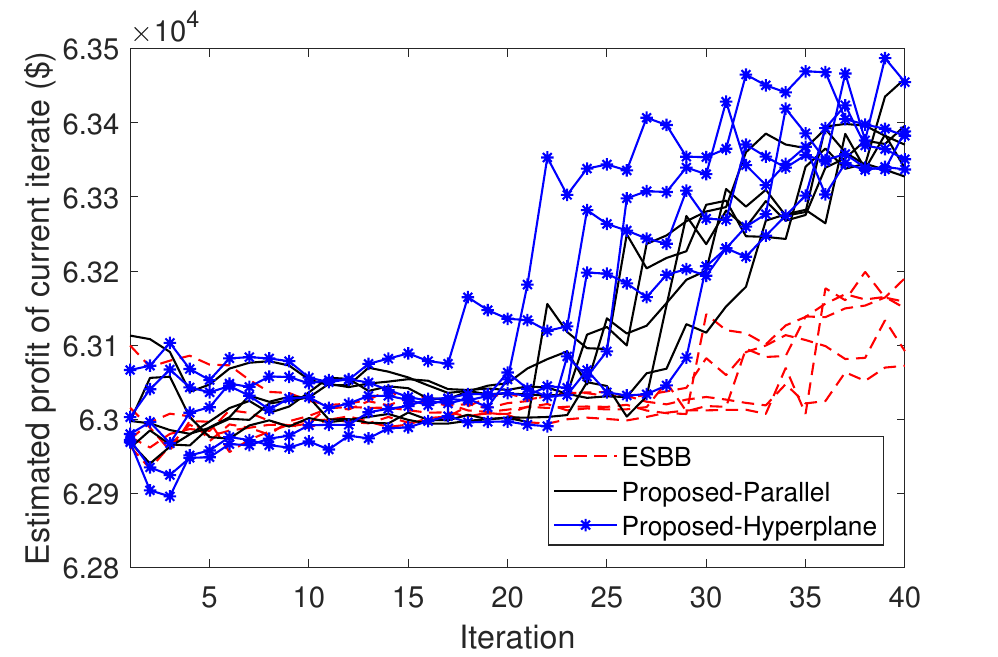}
		\caption{The objective function estimate of the current iterate across iterations for low-demand experiment.}
		\label{fig:zipcar_low}
		\endminipage\hfill
		\minipage{1\textwidth}
		\centering
		\includegraphics[width=0.9\linewidth]{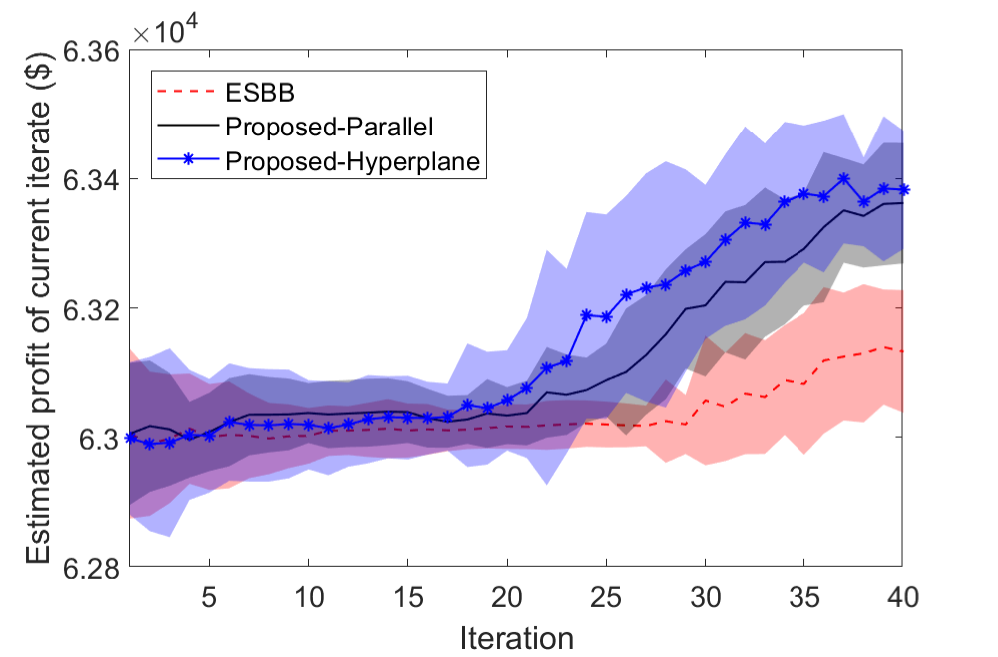}
		\caption{Objective function estimate of the current iterate across iterations(averaged over the 5 algorithm runs) for the low-demand experiment.}
		\label{fig:zipcar_low_agg}
		\endminipage\hfill
	\end{figure}

	Figure \ref{fig:zipcar_low} compares the objective estimate across iterations for each algorithm run of the low-demand experiment. 
	The x-axis displays the iteration index and the y-axis displays the performance estimate of the current iterate (i.e., simulation-based estimate of the objective function of the best point).
	The red dashed lines displays the results for the ESB\&B algorithm, the black solid lines for the proposed algorithm with parallel partition, and the blue asterisk lines for the proposed algorithm with hyperplane partition.
	Figure \ref{fig:zipcar_low_agg} displays the aggregated average performance for each algorithm. It plots the objective value of estimated optimal solution at each iteration (averaged over 5 algorithm runs) for the ESB\&B algorithm (red line with shaded $95\%$ confidence boundary), the proposed algorithm with parallel partition (black line with shaded $95\%$ confidence boundary), and the proposed algorithm with parallel partition (blue asterisk line with shaded $95\%$ confidence boundary).
	In both figures, we observe the following trends. 
	Initially, all algorithms start at a position similar to each other, which is the warm-start solution. 
	As iteration advances, both the proposed algorithms with parallel and hyperplane partition identify solutions with improved performance faster than the ESB\&B algorithm. 
	
	To further evaluate the performances of the solutions derived, we simulate each derived solution for $50$ replications.
	Table \ref{tab:low-demand} displays the performance statistics of the derived solutions from different algorithm runs. The first and second columns show the algorithm name and run id. The third and fourth columns show the mean and standard deviation of profits generated by each derived solution.
	We conduct a two-sample t-test between each derived solution from the proposed algorithm (with parallel and hyperplane partitions) and each derived solution from the ESB\&B algorithm. 
	The null hypothesis is that the average profit generated by the derived solutions from both algorithms are the same. 
	The corresponding alternative hypothesis is that the average profit generated by the solution derived from the proposed algorithm is higher.
	Table \ref{tab:low-demand-1} and \ref{tab:low-demand-2} display the p-values for the $50$ tests. The null hypothesis cannot be rejected if the p-value is greater than the significant level $0.05$ (colored red in the tables).
	The null hypothesis is rejected for all tests, except for the one which compares the derived solution from run $4$ of the ESB\&B algorithm and that from run $3$ of the proposed algorithm with parallel partition. 
	In other words, at the end of last iteration, the performances of the solutions derived from the proposed algorithm are mostly better than those derived from the ESB\&B algorithm.
	
	To investigate the added value of hyperplane cuts in the proposed algorithm, 
	we conduct a two-sample t-test for each solution derived from the proposed algorithm with parallel partition and each solution derived from the proposed algorithm with hyperplane partition.
	As before, our null hypothesis is that the average profits generated by both solutions are the same, and the alternative hypothesis is that the average profit generated by the solution derived from the proposed algorithm with hyperplane partition is higher.
	Table \ref{tab:low-demand-3} displays the p-values for all $25$ two-sample t-tests. The red cell in the table represents the test in which the null hypothesis cannot be rejected at significance level $0.05$ and there are 10 of those.
	The rest $15$ tests reject the null hypothesis and suggest that solution derived by the proposed algorithm with hyperplane partition generates higher profits. 

	\begin{table}[htp]
		\centering
		\begin{tabular}{c|ccc}
			\toprule
			Algorithm & Run & Mean & Standard deviation \\ \hline
			\multirow{5}{*}{ESB\&B algorithm} & 1 & 63066.95 & 173.73 \\
			& 2 & 63001.90 & 174.63 \\
			& 3 & 63056.80 & 175.06 \\
			& 4 & 63138.56 & 184.42 \\
			& 5 & 63076.53 & 174.82 \\ \hline
			\multirow{5}{*}{\begin{tabular}[c]{@{}c@{}}
					Proposed algorithm \\
					Parallel partition
			\end{tabular}} & 1 & 63204.83 & 175.60\\
			& 2 & 63220.92 & 161.94 \\
			& 3 & 63169.15 & 180.09  \\
			& 4 & 63244.88 & 150.86  \\
			& 5 & 63223.98 & 169.94 \\ \hline
			\multirow{5}{*}{\begin{tabular}[c]{@{}c@{}}
					Proposed algorithm  \\
					Hyperplane partition
			\end{tabular}} & 1 & 63260.53 & 204.98 \\
			& 2 & 63265.89 & 174.31 \\
			& 3 & 63247.42 & 168.43  \\
			& 4 & 63302.07 & 188.59 \\
			& 5 & 63322.54 & 180.93\\ \bottomrule
		\end{tabular}
		\caption{\label{tab:low-demand} Performance statistics of the derived final solutions in the low-demand case.} 	
	\end{table}
	
	\begin{table}[htp]
		\centering
		\begin{tabular}{c|c|ccccc}
			\toprule
			\multicolumn{1}{l|}{}   &     & \multicolumn{5}{c}{Proposed algorithm with parallel partition} \\ \hline
			\multicolumn{1}{l|}{}   & Run &   1    &   2    &            3            &   4    &     5     \\ \hline
			\multirow{5}{*}{ESB\&B} &  1  & 0.0001 & 0.0000 &         0.0024          & 0.0000 &  0.0000   \\
			&  2  & 0.0000 & 0.0000 &         0.0000          & 0.0000 &  0.0000   \\
			&  3  & 0.0000 & 0.0000 &         0.0010          & 0.0000 &  0.0000   \\
			&  4  & 0.0344 & 0.0098 & \textcolor{red}{0.2017} & 0.0011 &  0.0089   \\
			&  5  & 0.0002 & 0.0000 &         0.0052          & 0.0000 &  0.0000   \\ \bottomrule
		\end{tabular}
		\caption{\label{tab:low-demand-1} P-values of the two-sample t-test comparing the solutions derived by ESB\&B algorithm and the proposed algorithm with parallel partition in the low-demand experiment.}
	\end{table}
	
	\begin{table}[htp]
		\centering
		\begin{tabular}{c|c|ccccc}
			\toprule
			\multicolumn{1}{l|}{}  &     & \multicolumn{5}{c}{Proposed algorithm with hyperplane partition} \\ \hline
			\multicolumn{1}{l|}{}  & Run &   1    &   2    &   3    &   4    &              5               \\ \hline
			\multirow{5}{*}{ESB\&B} &  1  & 0.0000 & 0.0000 & 0.0000 & 0.0000 &            0.0000            \\
			&  2  & 0.0000 & 0.0000 & 0.0000 & 0.0000 &            0.0000            \\
			&  3  & 0.0000 & 0.0000 & 0.0000 & 0.0000 &            0.0000            \\
			&  4  & 0.0011 & 0.0003 & 0.0013 & 0.0000 &            0.0000            \\
			&  5  & 0.0000 & 0.0000 & 0.0000 & 0.0000 &            0.0000            \\ \bottomrule
		\end{tabular}
		\caption{\label{tab:low-demand-2} P-values of the two-sample t-test comparing the solutions derived by ESB\&B algorithm and the proposed algorithm with hyperplane partition in the low-demand experiment.}
	\end{table}
	
	\begin{table}[htp]
		\centering
		\begin{tabular}{c|c|ccccc}
			\toprule
			\multicolumn{1}{l|}{}  &     & \multicolumn{5}{c}{Proposed algorithm with hyperplane partition} \\ \hline
			\multicolumn{1}{l|}{}  & Run &   1    &   2    &   3    &   4    &              5               \\ \hline
			\multirow{5}{*}{\begin{tabular}[c]{@{}c@{}}
					Proposed algorithm \\
					with parallel partition
			\end{tabular}} &  1  & \textcolor{red}{0.0739} & 0.0420 & \textcolor{red}{0.1093} & 0.0045 &            0.0007            \\
			&  2  & \textcolor{red}{0.1433} & \textcolor{red}{0.0922} & \textcolor{red}{0.2123} & 0.0116 &            0.0019            \\
			&  3  & 0.0099 & 0.0037 & 0.0135 & 0.0002 &            0.0000            \\
			&  4  & \textcolor{red}{0.3325} & \textcolor{red}{0.2604} & \textcolor{red}{0.4684} & 0.0487 &            0.0109            \\
			&  5  & \textcolor{red}{0.1672} & 0.0000 & \textcolor{red}{0.2451} & 0.0160 &            0.0030            \\ \bottomrule
		\end{tabular}
		\caption{\label{tab:low-demand-3} P-values of the two-sample t-test comparing the solutions derived by the proposed algorithm with parallel and hyperplane partition in the low-demand experiment.}
	\end{table}
	
\begin{figure}[htp]
		\minipage{1\textwidth}
		\centering
		\includegraphics[width=0.9\linewidth]{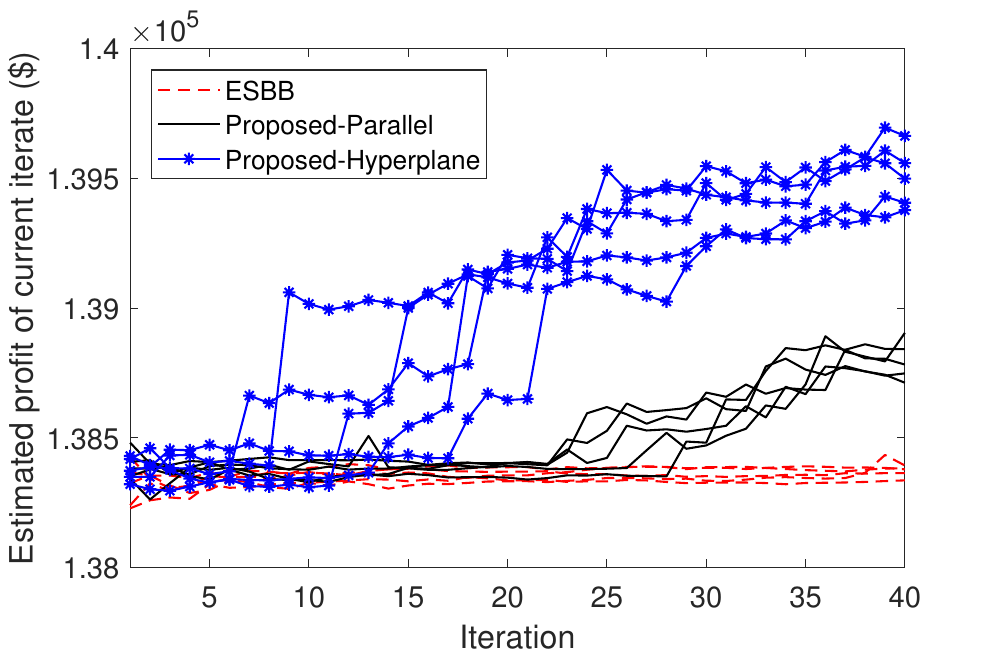}
		\caption{The objective function estimate of the current iterate across iterations for high-demand experiment.}
		\label{fig:zipcar_high}
		\endminipage\hfill
		\minipage{1\textwidth}
		\centering
		\includegraphics[width=0.9\linewidth]{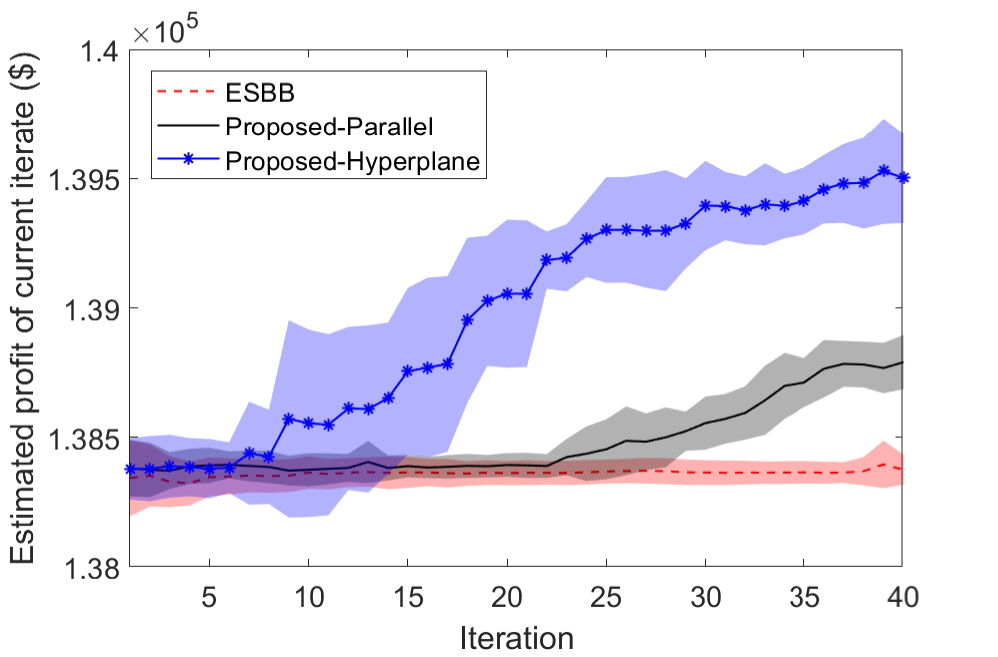}
		\caption{Objective function estimate of the current iterate across iterations (averaged over the 5 algorithm runs) for the high-demand experiment.}
		\label{fig:zipcar_high_agg}
		\endminipage\hfill
	\end{figure}

	Figure \ref{fig:zipcar_high} compares the objective estimate of the current iterate across iterations for each algorithm run of the high-demand experiment. This figure has a similar layout as Figure \ref{fig:zipcar_low}.
	Figure \ref{fig:zipcar_high_agg} displays the aggregated average performance for each algorithm. It has a similar layout as Figure \ref{fig:zipcar_low_agg}.
	From the beginning to the end, the ESB\&B algorithm does not find any solution with improved performance for four algorithm runs. 
	On the other hand, both the proposed algorithm with parallel and hyperplane partition do make improvement as iteration advances. 
	The proposed algorithm with hyperplane partition starts to find better solutions earlier than the one with parallel partition. 
	
	To further evaluate the performances of the solutions derived, we simulate each derived solution for $50$ replications.
	Table \ref{tab:high-demand} displays the performance statistics of the derived solutions from different algorithm runs. This table has a similar layout as Table \ref{tab:low-demand}.
	We first conduct a two-sample t-tests between each solution derived from the ESB\&B algorithm and each solution derived from the proposed algorithm (with parallel and hyperplane partition). 
	As before, the null hypothesis is that average profits generated by the solutions derived from both algorithms are the same, and the alternative hypothesis is that the solution derived from the proposed algorithm generates a higher profit. 
	Table \ref{tab:high-demand-1} and \ref{tab:high-demand-2} display the p-values for the $50$ tests. For all tests, the null hypothesis is rejected at significance level $0.05$. In other words, the proposed algorithm ends up with solutions with higher profits.
	
	To validate the added value of hyperplane cuts in the proposed algorithm, we conduct a two-sample t-test between each solution derived with parallel partition and each solution derived with hyperplane partition. 
	The null hypothesis is that the average profits generated by both solutions are the same. The alternative hypothesis is that the solution derived by the proposed algorithm with hyperplane partition generates a higher profit. 
	Table \ref{tab:high-demand-3} displays the p-values for the $25$ two-sample t-tests. For all the $25$ t-tests, we reject the null hypothesis at significance level $0.05$.
	Hence, the proposed algorithm with hyperplane partition ends up with derived solutions with better performances than those derived by the proposed algorithm with parallel partition. 
	This obvious improvement in the proposed algorithm's finite-time performance is because under high demand condition, customer spillover from one station to another happens more frequently and hence a hyperplane cut based on the total number of vehicles in a cluster of nearby stations forms a more efficient way of dividing the feasible region.
	
	\begin{table}[htp]
		\centering
		\begin{tabular}{c|ccc}
			\toprule
			Algorithm & Run & Mean & Standard deviation \\ \hline
			\multirow{5}{*}{ESB\&B algorithm} & 1 & 138420.78 & 177.05 \\
			& 2 & 138420.06 & 190.27 \\
			& 3 & 138390.60 & 158.24 \\
			& 4 & 138357.82 & 171.41 \\
			& 5 & 138357.76 & 210.50 \\ \hline
			\multirow{5}{*}{\begin{tabular}[c]{@{}c@{}}
					Proposed algorithm \\
					Parallel partition
			\end{tabular}} & 1 & 138669.72 & 184.94\\
			& 2 & 138703.72 & 166.58 \\
			& 3 & 138723.02 & 201.85  \\
			& 4 & 138751.44 & 204.59  \\
			& 5 & 138666.18 & 206.07 \\ \hline
			\multirow{5}{*}{\begin{tabular}[c]{@{}c@{}}
					Proposed algorithm  \\
					Hyperplane partition
			\end{tabular}} & 1 & 139502.98 & 196.23 \\
			& 2 & 139340.78 & 206.50 \\
			& 3 & 139188.80 & 212.03  \\
			& 4 & 139468.24 & 224.51 \\
			& 5 & 139272.32 & 195.04\\ \bottomrule
		\end{tabular}
		\caption{\label{tab:high-demand} Performance statistics of the derived final solutions in the high-demand case.} 	
	\end{table}
	
	\begin{table}[htp]
		\centering
		\begin{tabular}{c|c|ccccc}
			\toprule
			\multicolumn{1}{l|}{}  &     & \multicolumn{5}{c}{Proposed algorithm with parallel partition} \\ \hline
			\multicolumn{1}{l|}{}  & Run &   1    &   2    &   3    &   4    &             5              \\ \hline
			\multirow{5}{*}{ESB\&B} &  1  & 0.0000 & 0.0000 & 0.0000 & 0.0000 &           0.0000           \\
			&  2  & 0.0000 & 0.0000 & 0.0000 & 0.0000 &           0.0000           \\
			&  3  & 0.0000 & 0.0000 & 0.0000 & 0.0000 &           0.0000           \\
			&  4  & 0.0000 & 0.0000 & 0.0000 & 0.0000 &           0.0000           \\
			&  5  & 0.0000 & 0.0000 & 0.0000 & 0.0000 &           0.0000           \\ \bottomrule
		\end{tabular}
		\caption{\label{tab:high-demand-1} P-values of the two-sample t-test comparing the solutions derived by ESB\&B algorithm and the proposed algorithm with parallel partition in the high-demand experiment.}
	\end{table}
	
	\begin{table}[htp]
		\centering
		\begin{tabular}{c|c|ccccc}
			\toprule
			\multicolumn{1}{l|}{}  &     & \multicolumn{5}{c}{Proposed algorithm with hyperplane partition} \\ \hline
			\multicolumn{1}{l|}{}  & Run &   1    &   2    &   3    &   4    &              5               \\ \hline
			\multirow{5}{*}{ESB\&B} &  1  & 0.0000 & 0.0000 & 0.0000 & 0.0000 &            0.0000            \\
			&  2  & 0.0000 & 0.0000 & 0.0000 & 0.0000 &            0.0000            \\
			&  3  & 0.0000 & 0.0000 & 0.0000 & 0.0000 &            0.0000            \\
			&  4  & 0.0000 & 0.0000 & 0.0000 & 0.0000 &            0.0000            \\
			&  5  & 0.0000 & 0.0000 & 0.0000 & 0.0000 &            0.0000            \\ \bottomrule
		\end{tabular}
		\caption{\label{tab:high-demand-2} P-values of the two-sample t-test comparing the solutions derived by ESB\&B algorithm and the proposed algorithm with hyperplane partition in the high-demand experiment.}
	\end{table}
	
	\begin{table}[htp]
		\centering
		\begin{tabular}{c|c|ccccc}
			\toprule
			\multicolumn{1}{l|}{}  &     & \multicolumn{5}{c}{Proposed algorithm with hyperplane partition} \\ \hline
			\multicolumn{1}{l|}{}  & Run &   1    &   2    &   3    &   4    &              5               \\ \hline
			\multirow{5}{*}{\begin{tabular}[c]{@{}c@{}}
					Proposed algorithm \\
					with parallel partition
			\end{tabular}} &  1  & 0.0000 & 0.0000 & 0.0000 & 0.0000 &            0.0000            \\
			&  2  & 0.0000 & 0.0000 & 0.0000 & 0.0000 &            0.0000            \\
			&  3  & 0.0000 & 0.0000 & 0.0000 & 0.0000 &            0.0000            \\
			&  4  & 0.0000 & 0.0000 & 0.0000 & 0.0000 &            0.0000            \\
			&  5  & 0.0000 & 0.0000 & 0.0000 & 0.0000 &            0.0000            \\ \bottomrule
		\end{tabular}
		\caption{\label{tab:high-demand-3} P-values of the two-sample t-test comparing the solutions derived by the proposed algorithm with parallel and hyperplane partition in the high-demand experiment.}
	\end{table}

	In summary, compared to the ESB\&B algorithm, the proposed algorithm can improve the finite-time performances by exploring the underlying function structure through sampled points and incorporating prior knowledge of the problem-specific structures.
	
	\section{Conclusion}\label{sec:conclude}
	In this paper, we propose an adaptive partitioning strategy and combine it with the ESB\&B framework developed by \citet{xu2013empirical}. This proposed partitioning strategy is a sample-driven approach that explores the structure of the underlying objective function by iteratively dividing the feasible region into subregions such that sampled points within the same subregion have similar performances. 
	The proposed partitioning strategy can be integrated with prior knowledge of specific problem structures to form more efficient hyperplane partitions.  
	Solving the proposed partitioning problem is fast and efficient through solving the formulated MIP via the local search algorithm developed by \citet{dunn2018optimal}. 
	The proposed algorithm combines the proposed adaptive partitioning strategy within the ESB\&B framework. 
	It is a general-purpose algorithm that converges globally for discrete SO problems with finite and convex feasible region, and it can be applied to general SO problems without significant amount of modification.
	The proposed algorithm improves the finite-time performances of the original ESB\&B algorithm with equal partitioning strategy.
	
	For SO problems with tight sampling budget at each iteration, a smart sampling strategy other than uniform sampling is important. This can be a potential research direction to explore. 
	Currently, the adaptive partitioning strategy supports parallel partition and only hyperplane partitions with potential hyperplane cuts known in advance.
	It is an interesting research direction to automatically generate hyperplane cuts for any problems without such prior knowledge or such prior knowledge indicates non-hyperplane relations \citep{lu2018probabilistic, lu2022analytical, lu2024link}.

\bibliographystyle{elsarticle-harv}
\bibliography{SObiblio}	
\end{document}